\documentclass[12pt]{article}%
\usepackage{amsmath}
\usepackage{amsfonts}
\usepackage{amssymb}
\usepackage{graphicx}%
\providecommand{\U}[1]{\protect\rule{.1in}{.1in}}
\newtheorem{theorem}{Theorem}[section]

\newtheorem{proposition}[theorem]{Proposition}
\newtheorem{remark}[theorem]{Remark}

\numberwithin{equation}{section}

\begin{document}

\title{Approximate direct and inverse scattering for the AKNS system}
\author{Vladislav V. Kravchenko\\{\small Departamento de Matem\'{a}ticas, Cinvestav, Unidad Quer\'{e}taro, }\\{\small Libramiento Norponiente \#2000, Fracc. Real de Juriquilla,
Quer\'{e}taro, Qro., 76230 Mexico,}\\{\small Regional Mathematical Center, Southern Federal University,
Rostov-on-Don, 344090, Russia}\\{\small e-mail: vkravchenko@math.cinvestav.edu.mx}}
\maketitle

\begin{abstract}
We study the direct and inverse scattering problems for the AKNS
(Ablowitz-Kaup-Newell-Segur) system. New representations for the Jost
solutions are obtained in the form of the power series in terms of a
transformed spectral parameter. In terms of that parameter, the Jost solutions
are convergent power series in corresponding unit disks. For the coefficients
of the series simple recurrent integration procedures are devised.

Solution of the direct scattering problem reduces to computing the
coefficients and locating zeros of corresponding analytic functions in the
interior of the unit disk.

Solution of the inverse scattering problem reduces to the solution of two
systems of linear algebraic equations for the power series coefficients, while
the potentials are recovered from the first coefficients.

The overall approach leads to a simple and efficient method for the numerical
solution of both direct and inverse scattering problems, which is illustrated
by numerical examples.\ 

\end{abstract}

\section{Introduction}

We consider the direct scattering and inverse scattering problems for the AKNS
(Ablowitz-Kaup-Newell-Segur) system \cite{AKNS}
\begin{align*}
\frac{dn_{1}(x)}{dx}+i\rho n_{1}(x) &  =q(x)n_{2}(x),\\
\frac{dn_{2}(x)}{dx}-i\rho n_{2}(x) &  =r(x)n_{1}(x),\quad x\in(-\infty
,\infty),
\end{align*}
where $q(x)$ and $r(x)$ are complex valued functions, integrable and of a
sufficiently fast decay at infinity, $\rho\in\mathbb{C}$ is a spectral
parameter. This includes the Zakharov-Shabat system \cite{ZS} by restricting
to the case where $r(x)=\pm\overline{q(x)}$, as well as the system considered
by Gardner-Green-Kruskal-Miura \cite{Gardner et al} by setting $q\equiv-1$.
Main results of the scattering theory for the AKNS system are well known and
presented in several books (see, e.g., \cite{Ablowitz Segur}, \cite{Gerdjikov
et al 2008}, \cite{Lamb}). Direct and inverse scattering problems for the AKNS
system encounter multiple applications. They are involved in solving nonlinear
evolution equations, such as the nonlinear Schr\"{o}dinger, sine-Gordon,
Korteweg - de Vries, modified Korteweg - de Vries equations among others
\cite{Ablowitz Clarkson}, \cite{Ablowitz Segur}, \cite{Beals et al 1988},
\cite{Deift et al 1979}, \cite{Drazin et al 1989}, \cite{Lamb}, and also for
modelling other physical phenomena (see, e.g., \cite{Sacks}, \cite{Tang et al}).

Numerical methods for solving both direct and inverse scattering problems for
the AKNS or Zakharov-Shabat systems are subject of numerous publications (see,
e.g., \cite{Arico et al 2011}, \cite{Chimmalgi et al 2019}, \cite{Delitsyn
2022}, \cite{Fermo et al 2016}, \cite{Frangos et al 1991}, \cite{Frumin et al
2015}, \cite{Gorbenko et al}, \cite{Le et al 2014}, \cite{Medvedev et al
2023}, \cite{Mullyadzhanov et al 2021}, \cite{Trogdon2021}, \cite{Trogdon
Olver}, \cite{Trogdon et al 2012}, \cite{Turitsyn et al 2017}, \cite{Wahls et
al 2018}, \cite{Xiao et al 2002}). However, the problems are computationally
challenging, and purely numerical techniques reveal limitations both in
accuracy and speed.

In the present work an entirely different approach to solving scattering and
inverse scattering problems for the AKNS system is developed. Let us outline
its important features. The AKNS system possesses four fundamental solutions
called Jost solutions, which satisfy corresponding asymptotic relations at
infinity. For example, one of the Jost solutions satisfies the asymptotics
\[
\psi(\rho,x)=\left(
\begin{tabular}
[c]{c}%
$\psi_{1}(\rho,x)$\\
$\psi_{2}(\rho,x)$%
\end{tabular}
\right)  \sim\left(
\begin{tabular}
[c]{c}%
$0$\\
$1$%
\end{tabular}
\ \right)  e^{i\rho x},\quad x\rightarrow\infty.
\]

It is known that it admits the integral representation
\begin{equation}
\psi(\rho,x)=\left(
\begin{tabular}
[c]{c}%
$0$\\
$1$%
\end{tabular}
\ \right)  e^{i\rho x}+\int_{x}^{\infty}A(x,t)e^{i\rho t}dt,
\label{psi integral}%
\end{equation}
where the vector kernel $A(x,\cdot)$ belongs to $\mathcal{L}_{2}\left(
x,\infty\right)  $ componentwise. Analogous results are available for the rest
of the Jost solutions. As a consequence of the square integrability of
$A(x,\cdot)$, it can be expanded into the Fourier-Laguerre series%
\begin{equation}
A(x,t)=\sum_{n=0}^{\infty}a_{n}(x)L_{n}(t-x)e^{\frac{x-t}{2}},
\label{A series}%
\end{equation}
where $L_{n}$ stands for the Laguerre polynomial of order $n$, and each
coefficient is a vector function,
\[
a_{n}(x)=\left(
\begin{tabular}
[c]{c}%
$a_{1,n}(x)$\\
$a_{2,n}(x)$%
\end{tabular}
\right)  .
\]
Substitution of (\ref{A series}) into (\ref{psi integral}) leads to the series
representation for the Jost solution
\[
\psi(\rho,x)=e^{i\rho x}\left(  \left(
\begin{tabular}
[c]{c}%
$0$\\
$1$%
\end{tabular}
\ \right)  +\left(  z+1\right)  \sum_{n=0}^{\infty}\left(  -1\right)
^{n}z^{n}a_{n}(x)\right)  ,
\]
where $z$ arises in a natural way as a result of integration of the products
of the Laguerre polynomials with the exponential functions. It is defined by
\[
z=\frac{\frac{1}{2}+i\rho}{\frac{1}{2}-i\rho},
\]
which is nothing but a M\"{o}bius map transforming the upper half-plane
$\operatorname{Im}\rho\geq0$ onto the unit disk $\left\vert z\right\vert
\leq1$. As a function of $\rho$ the Jost solution $\psi(\rho,x)$ is analytic
in the upper half-plane, so in terms of $z$ it is analytic in the unit disk.
Now, the crucial question is whether it is possible to efficiently compute the
coefficients $a_{n}(x)$. The answer is positive. We develop a recurrent
integration procedure for computing $a_{n}(x)$. It requires to obtain the Jost
solution $\psi(\frac{i}{2},x)$ in the first step, which gives us the first
coefficient
\begin{equation}
a_{0}(x)=e^{\frac{x}{2}}\psi(\frac{i}{2},x)-\left(
\begin{tabular}
[c]{c}%
$0$\\
$1$%
\end{tabular}
\right)  . \label{a0=}%
\end{equation}
The subsequent coefficients are computed by the procedure which involves only
multiplication by known functions and integration. The procedure is simple and
allows us to compute hundreds of the coefficients $a_{n}(x)$, if necessary.
Analogous representations and recurrent procedures are obtained for the other
Jost solutions.

Thus, the Jost solutions are represented in the form of power series with
respect to the transformed spectral parameter. We refer to these power series
expansions of the Jost solutions as the spectral parameter power series
(SPPS), and to their coefficients as the SPPS coefficients.

Therefore, after having computed the SPPS coefficients for $x=0$, the
computation of the scattering data, such as the entries of the scattering
(transfer) matrix for $\rho\in\mathbb{R}$ reduces to the polynomial evaluation
along the unit circle ($\left\vert z\right\vert =1$ corresponds to $\rho
\in\mathbb{R}$), while the computation of the eigenvalues reduces to the
location of zeros of corresponding polynomials inside the unit disk. We show
that all this works numerically and gives accurate results outdoing best
available techniques.

The obtained SPPS representations lend themselves to an elementary approach to
solving the inverse scattering problem. Here another feature of the SPPS
coefficients results to be crucial. The first coefficient of any of the four
SPPS representations (for the four Jost solutions) is sufficient for
recovering both potentials $q(x)$ and $r(x)$. This is obvious, for example,
from (\ref{a0=}). Indeed, if $a_{0}(x)$ is known, then from (\ref{a0=}) we
obtain $\psi(\frac{i}{2},x)$, and hence $q(x)$ and $r(x)$ can be computed by
using the AKNS system as
\[
q(x)=\frac{\psi_{1}^{\prime}(\frac{i}{2},x)-\psi_{1}(\frac{i}{2},x)/2}%
{\psi_{2}(\frac{i}{2},x)},\quad r(x)=\frac{\psi_{2}^{\prime}(\frac{i}%
{2},x)+\psi_{2}(\frac{i}{2},x)/2}{\psi_{1}(\frac{i}{2},x)}.
\]
To compute the SPPS coefficients we construct two systems of linear algebraic
equations directly from the scattering relations between the Jost solutions,
one system corresponds to the first components of the Jost solutions, while
the other to their second components.

Thus, the implementation of our approach to the inverse scattering problem
consists in substituting the SPPS representations for the Jost solutions into
the scattering relations and writing down the arising two systems of linear
algebraic equations for the SPPS coefficients. Solving the systems we recover
the potentials from the very first SPPS coefficients.

The idea of using series representations for the integral kernels of the Jost
solutions was proposed first in \cite{Kr2019MMAS InverseAxis} in the context
of the Schr\"{o}dinger equation and was developed further in subsequent
publications \cite{DKK2019MMAS}, \cite{DKK Jost}, \cite{GKT2023MMAS},
\cite{KarapetyantsKravchenkoBook}, \cite{KrBook2020}. For the Zakharov-Shabat
system with a real-valued potential it was developed in \cite{Kr2025ZSreal},
where the reduction of the system to a couple of Schr\"{o}dinger equations was
used, so that the results from the preceding work \cite{Kr2019MMAS
InverseAxis}, \cite{GKT2023MMAS} could be adapted. The general AKNS system
considered in the present work is quite more challenging. Its consideration
required the whole construction to be developed from the ground up.

Thus, the novelty of the present work is a new series representation for the
Jost solutions of the AKNS system and based on it a new approach to numerical
solution of direct and inverse scattering problems. The functionality of the
approach is illustrated by several numerical tests.

In Section 2 we introduce the necessary notations and formulate the problems.
In Section 3 we introduce the SPPS representations for the Jost solutions and
prove their main properties. In Section 4 we develop the recurrent integration
procedure for computing the SPPS coefficients. Sections 5 and 6 are dedicated
to solving the direct and inverse scattering problems, respectively. In
Section 7 we discuss the numerical implementation of the approach and
illustrate it with several numerical examples. Finally, Section 8 contains
some concluding remarks.

\section{Direct and inverse problems, preliminaries}

Let $q(x)$ and $r(x)$ be complex valued, integrable functions defined on
$\left(  -\infty,\infty\right)  $. Consider the AKNS (or generalized
Zakharov-Shabat) system
\begin{equation}
\frac{dn_{1}(x)}{dx}+i\rho n_{1}(x)=q(x)n_{2}(x),\label{ZS1}%
\end{equation}%
\begin{equation}
\frac{dn_{2}(x)}{dx}-i\rho n_{2}(x)=r(x)n_{1}(x)\label{ZS2}%
\end{equation}
where $\rho\in\mathbb{C}$ is a spectral parameter.

For any two solutions of (\ref{ZS1}), (\ref{ZS2})
\[
\theta(x)=\left(
\begin{tabular}
[c]{c}%
$\theta_{1}(x)$\\
$\theta_{2}(x)$%
\end{tabular}
\right)  \quad\text{and}\quad\omega(x)=\left(
\begin{tabular}
[c]{c}%
$\omega_{1}(x)$\\
$\omega_{2}(x)$%
\end{tabular}
\right)  ,
\]
for the same value of the parameter $\rho$, the expression $W\left[
\theta;\omega\right]  :=\theta_{1}\omega_{2}-\theta_{2}\omega_{1}$ is
constant, $\frac{dW}{dx}=0$, and the nonvanishing of $W\left[  \theta
;\omega\right]  $ guarantees the linear independence of $\theta$ and $\omega$.

It is convenient to consider solutions of (\ref{ZS1}), (\ref{ZS2}) as
functions of $\rho$ as well, so we write
\[
n=n(\rho,x)=\left(
\begin{tabular}
[c]{c}%
$n_{1}(\rho,x)$\\
$n_{2}(\rho,x)$%
\end{tabular}
\right)  .
\]

We will suppose that $q(x)$ and $r(x)$ decay sufficiently fast at $\pm\infty$,
so that there exist the unique, so-called Jost solutions
\[
\varphi(\rho,x)\sim\left(
\begin{tabular}
[c]{c}%
$1$\\
$0$%
\end{tabular}
\ \right)  e^{-i\rho x},\quad\widetilde{\varphi}(\rho,x)=\left(
\begin{tabular}
[c]{c}%
$0$\\
$-1$%
\end{tabular}
\ \right)  e^{i\rho x},\quad x\rightarrow-\infty
\]
and%
\[
\psi(\rho,x)\sim\left(
\begin{tabular}
[c]{c}%
$0$\\
$1$%
\end{tabular}
\ \right)  e^{i\rho x},\quad\widetilde{\psi}(\rho,x)\sim\left(
\begin{tabular}
[c]{c}%
$1$\\
$0$%
\end{tabular}
\ \right)  e^{-i\rho x},\quad x\rightarrow\infty,
\]
and the integral representations of the form (\ref{phi Levin}%
)-(\ref{psitil Levin}) below are valid.

We have
\begin{equation}
W\left[  \varphi(\rho,x);\widetilde{\varphi}(\rho,x)\right]  =-1\quad
\text{and\quad}W\left[  \psi(\rho,x);\widetilde{\psi}(\rho,x)\right]  =-1.
\label{Wronski}%
\end{equation}
Thus, solutions $\varphi(\rho,x)$ and $\widetilde{\varphi}(\rho,x)$ as well as
$\psi(\rho,x)$ and $\widetilde{\psi}(\rho,x)$ are linearly independent, and
there exist the scalars $\mathbf{a}(\rho)$, $\widetilde{\mathbf{a}}(\rho)$,
$\mathbf{b}(\rho)$, $\widetilde{\mathbf{b}}(\rho)$, such that
\begin{equation}
\varphi(\rho,x)=\mathbf{a}(\rho)\widetilde{\psi}(\rho,x)+\mathbf{b}(\rho
)\psi(\rho,x), \label{rel 1}%
\end{equation}%
\begin{equation}
\widetilde{\varphi}(\rho,x)=-\widetilde{\mathbf{a}}(\rho)\psi(\rho
,x)+\widetilde{\mathbf{b}}(\rho)\widetilde{\psi}(\rho,x). \label{rel 2}%
\end{equation}

Here we follow the notation from \cite{Ablowitz Segur}. The scattering
(transfer) matrix is defined as
\[
S(\rho)=\left(
\begin{tabular}
[c]{cc}%
$\mathbf{a}(\rho)$ & $\mathbf{b}(\rho)$\\
$\widetilde{\mathbf{b}}(\rho)$ & $-\widetilde{\mathbf{a}}(\rho)$%
\end{tabular}
\right)  .
\]
The relation holds%
\begin{equation}
\mathbf{a}(\rho)\widetilde{\mathbf{a}}(\rho)+\mathbf{b}(\rho)\widetilde
{\mathbf{b}}(\rho)=1. \label{aatil}%
\end{equation}

In terms of the Jost solutions the entries of the scattering matrix can be
written as follows \cite{Ablowitz Segur}%
\begin{equation}
\mathbf{a}(\rho)=W\left[  \varphi(\rho,0);\psi(\rho,0)\right]  ,\qquad
\widetilde{\mathbf{a}}(\rho)=W\left[  \widetilde{\varphi}(\rho,0);\widetilde
{\psi}(\rho,0)\right]  , \label{atila}%
\end{equation}%
\begin{equation}
\mathbf{b}(\rho)=-W\left[  \varphi(\rho,0);\widetilde{\psi}(\rho,0)\right]
,\qquad\widetilde{\mathbf{b}}(\rho)=W\left[  \widetilde{\varphi}(\rho
,0);\psi(\rho,0)\right]  . \label{btilb}%
\end{equation}

System (\ref{ZS1}), (\ref{ZS2}) may admit a number of the eigenvalues, which
coincide with zeros of $\mathbf{a}(\rho)$ in $\mathbb{C}^{+}:=\left\{  \rho
\in\mathbb{C}\mid\operatorname{Im}\rho>0\right\}  $ and zeros of
$\widetilde{\mathbf{a}}(\rho)$ in $\mathbb{C}^{-}:=\left\{  \rho\in
\mathbb{C}\mid\operatorname{Im}\rho<0\right\}  $. We will assume that the
number of the eigenvalues is finite. If $\mathbf{a}(\rho_{m})=0$,
$\operatorname{Im}\rho_{m}>0$, then $\varphi(\rho,x)$ and $\psi(\rho,x)$ are
linearly dependent, hence there exists a multiplier constant $c_{m}\neq0$ such
that
\begin{equation}
\varphi(\rho_{m},x)=c_{m}\psi(\rho_{m},x).\label{eig1}%
\end{equation}
$c_{m}$ is called the norming constant associated with the eigenvalue
$\rho_{m}$. Analogously, if $\widetilde{\mathbf{a}}(\widetilde{\rho}_{m})=0$,
$\operatorname{Im}\widetilde{\rho}_{m}<0$, there exists the norming constant
$\widetilde{c}_{m}\neq0$ such that%
\begin{equation}
\widetilde{\varphi}(\widetilde{\rho}_{m},x)=\widetilde{c}_{m}\widetilde{\psi
}(\widetilde{\rho}_{m},x).\label{eig2}%
\end{equation}
The number of the eigenvalues in $\mathbb{C}^{+}$ we denote by $M$, while the
number of the eigenvalues in $\mathbb{C}^{-}$ we denote by $\widetilde{M}$.

The direct scattering problem consists in finding $S(\rho)$ for all $\rho
\in\mathbb{R}$ as well as in finding all the eigenvalues and norming
constants. In other words, when solving the direct scattering problem for
system (\ref{ZS1}), (\ref{ZS2}) one needs to find the set of the scattering
data%
\begin{equation}
SD:=\left\{  S(\rho),\,\rho\in\mathbb{R};\left\{  \rho_{m},c_{m}\right\}
_{m=1}^{M},\left\{  \widetilde{\rho}_{m},\widetilde{c}_{m}\right\}
_{m=1}^{\widetilde{M}}\right\}  . \label{SD}%
\end{equation}

The inverse scattering problem consists in recovering the coefficients $q(x)$
and $r(x)$ in (\ref{ZS1}), (\ref{ZS2}) for all $x\in\mathbb{R}$ from the set
of the scattering data $SD$.

\section{Series representations for Jost solutions}

\subsection{Integral representations, series expansions of integral kernels}

As we stated earlier, we assume that $q(x)$, $r(x)$ decay sufficiently fast,
when $x\rightarrow\pm\infty$, that the Jost solutions admit the integral
representations%
\begin{equation}
\varphi(\rho,x)=\left(
\begin{tabular}
[c]{c}%
$1$\\
$0$%
\end{tabular}
\ \right)  e^{-i\rho x}+\int_{-\infty}^{x}B(x,t)e^{-i\rho t}dt,
\label{phi Levin}%
\end{equation}%
\begin{equation}
\widetilde{\varphi}(\rho,x)=\left(
\begin{tabular}
[c]{c}%
$0$\\
$-1$%
\end{tabular}
\ \right)  e^{i\rho x}+\int_{-\infty}^{x}\widetilde{B}(x,t)e^{i\rho t}dt,
\label{phitil Levin}%
\end{equation}%
\begin{equation}
\psi(\rho,x)=\left(
\begin{tabular}
[c]{c}%
$0$\\
$1$%
\end{tabular}
\ \right)  e^{i\rho x}+\int_{x}^{\infty}A(x,t)e^{i\rho t}dt, \label{psi Levin}%
\end{equation}%
\begin{equation}
\widetilde{\psi}(\rho,x)=\left(
\begin{tabular}
[c]{c}%
$1$\\
$0$%
\end{tabular}
\ \right)  e^{-i\rho x}+\int_{x}^{\infty}\widetilde{A}(x,t)e^{-i\rho t}dt,
\label{psitil Levin}%
\end{equation}
where the vector kernels $A(x,\cdot)$, $\widetilde{A}(x,\cdot)$ belong to
$\mathcal{L}_{2}\left(  x,\infty\right)  $ componentwise, and $B(x,\cdot)$,
$\widetilde{B}(x,\cdot)$ to $\mathcal{L}_{2}\left(  -\infty,x\right)  $. In
this case we say that the complex-valued integrable functions $q(x)$, $r(x)$,
tending to zero at $\pm\infty$, belong to the class $\mathcal{P}$. The exact
conditions for $q(x)$ and $r(x)$ to belong to $\mathcal{P}$ to our best
knowledge are unknown, however, in the important case $r(x)=-\overline{q}(x)$
it is sufficient to require $q(x)\in\mathcal{L}_{1}\left(  -\infty
,\infty\right)  \cap\mathcal{L}_{2}\left(  -\infty,\infty\right)  $
\cite{Hryniv Manko 2016}.

\begin{theorem}
Let $q(x)$ and $r(x)$ belong to the class $\mathcal{P}$. Then the integral
kernels in (\ref{phi Levin})-(\ref{psitil Levin}) admit the following series
representations
\begin{equation}
A(x,t)=\sum_{n=0}^{\infty}a_{n}(x)L_{n}(t-x)e^{\frac{x-t}{2}},\quad
\widetilde{A}(x,t)=\sum_{n=0}^{\infty}\widetilde{a}_{n}(x)L_{n}(t-x)e^{\frac
{x-t}{2}},\label{A}%
\end{equation}%
\begin{equation}
B(x,t)=\sum_{n=0}^{\infty}b_{n}(x)L_{n}(x-t)e^{-\frac{x-t}{2}},\quad
\widetilde{B}(x,t)=\sum_{n=0}^{\infty}\widetilde{b}_{n}(x)L_{n}(x-t)e^{-\frac
{x-t}{2}},\label{B}%
\end{equation}
where $L_{n}$ stands for the Laguerre polynomial of order $n$. Each
coefficient is a vector function, for example,
\[
a_{n}(x)=\left(
\begin{tabular}
[c]{c}%
$a_{1,n}(x)$\\
$a_{2,n}(x)$%
\end{tabular}
\ \ \right)  .
\]
For any $x$ fixed the series converge in mean square.
\end{theorem}

\textbf{Proof. }The proof of (\ref{A}) and (\ref{B}) is quite elementary and
similar to that from \cite{Kr2019MMAS InverseAxis}. Indeed, consider, for
example, (\ref{psi Levin}) and denote $a_{1}(x,t):=e^{\frac{t}{2}}%
A_{1}(x,x+t)$, where $A_{1}$ is the first component of the vector kernel $A$.
Since $A_{1}(x,\cdot)\in\mathcal{L}_{2}\left(  x,\infty\right)  $ we have that
$a_{1}(x,\cdot)$ belongs to $\mathcal{L}_{2}\left(  0,\infty;e^{-t}\right)  $.
Indeed,
\[
\int_{0}^{\infty}e^{-t}\left\vert a_{1}\right\vert ^{2}(x,t)dt=\int
_{0}^{\infty}\left\vert A_{1}\right\vert ^{2}(x,x+t)dt=\int_{x}^{\infty
}\left\vert A_{1}\right\vert ^{2}(x,t)dt<\infty.
\]
Hence $a_{1}(x,t)$ admits the representation (see, e.g., \cite{Olver et al
2010}, \cite{Suetin})
\[
a_{1}(x,t)=\sum_{n=0}^{\infty}a_{1,n}(x)L_{n}(t).
\]
For every $x\geq0$ the series converges in the norm of the space $L_{2}\left(
0,\infty;e^{-t}\right)  $. Returning to $A_{1}(x,t)$, we obtain the equality
\begin{equation}
A_{1}(x,t)=\sum_{n=0}^{\infty}a_{1,n}(x)L_{n}(t-x)e^{\frac{x-t}{2}}.\label{A1}%
\end{equation}
The proof for the second component of $A(x,t)$ is analogous, as well as for
the other three kernels in (\ref{A}) and (\ref{B}).$\qquad\square$

\subsection{Spectral parameter power series for Jost solutions}

\begin{theorem}
Let $q(x)$, $r(x)$ belong to the class $\mathcal{P}$. Then the corresponding
Jost solutions admit the series representations%
\begin{equation}
\varphi(\rho,x)=e^{-i\rho x}\left(  \left(
\begin{tabular}
[c]{c}%
$1$\\
$0$%
\end{tabular}
\right)  +\left(  z+1\right)  \sum_{n=0}^{\infty}\left(  -1\right)  ^{n}%
z^{n}b_{n}(x)\right)  , \label{phi=}%
\end{equation}%
\begin{equation}
\widetilde{\varphi}(\rho,x)=e^{i\rho x}\left(  \left(
\begin{tabular}
[c]{c}%
$0$\\
$-1$%
\end{tabular}
\right)  +\left(  \widetilde{z}+1\right)  \sum_{n=0}^{\infty}\left(
-1\right)  ^{n}\widetilde{z}^{n}\widetilde{b}_{n}(x)\right)  , \label{phitil=}%
\end{equation}%
\begin{equation}
\psi(\rho,x)=e^{i\rho x}\left(  \left(
\begin{tabular}
[c]{c}%
$0$\\
$1$%
\end{tabular}
\right)  +\left(  z+1\right)  \sum_{n=0}^{\infty}\left(  -1\right)  ^{n}%
z^{n}a_{n}(x)\right)  , \label{psi=}%
\end{equation}%
\begin{equation}
\widetilde{\psi}(\rho,x)=e^{-i\rho x}\left(  \left(
\begin{tabular}
[c]{c}%
$1$\\
$0$%
\end{tabular}
\right)  +\left(  \widetilde{z}+1\right)  \sum_{n=0}^{\infty}\left(
-1\right)  ^{n}\widetilde{z}^{n}\widetilde{a}_{n}(x)\right)  , \label{psitil=}%
\end{equation}
where the coefficients $a_{n}$, $\widetilde{a}_{n}$, $b_{n}$ and
$\widetilde{b}_{n}$ are those from (\ref{A}) and (\ref{B}), and
\begin{equation}
z=z\left(  \rho\right)  :=\frac{\frac{1}{2}+i\rho}{\frac{1}{2}-i\rho}%
,\quad\widetilde{z}=\widetilde{z}\left(  \rho\right)  :=\frac{\frac{1}%
{2}-i\rho}{\frac{1}{2}+i\rho}. \label{z}%
\end{equation}
The power series in (\ref{phi=}), (\ref{psi=}) converge in the unit disk
$D=\left\{  z\in\mathbb{C}:\,\left\vert z\right\vert <1\right\}  $ of the
complex plane of the variable $z$, and the power series in (\ref{phitil=}),
(\ref{psitil=}) converge in the unit disk $\widetilde{D}=\left\{
\widetilde{z}\in\mathbb{C}:\,\left\vert \widetilde{z}\right\vert <1\right\}  $
of the complex plane of the variable $\widetilde{z}$.
\end{theorem}

\textbf{Proof. }Substitution of expressions (\ref{A}), (\ref{B}) into
respective integral representations for the Jost solutions lead to the series
representations for the latter ones. Namely, let us prove, for example,
(\ref{psitil=}). Substitution of the series expansion of $\widetilde{A}(x,t)$
from (\ref{A}) into (\ref{psitil Levin}) leads to the equality%
\[
\widetilde{\psi}(\rho,x)=e^{-i\rho x}\left(  \left(
\begin{tabular}
[c]{c}%
$1$\\
$0$%
\end{tabular}
\ \right)  +\sum_{n=0}^{\infty}\widetilde{a}_{n}(x)\int_{0}^{\infty}%
L_{n}(t)e^{-\left(  \frac{1}{2}+i\rho\right)  t}dt\right)  ,
\]
where the change of the order of integration and summation is due to
Parseval's identity \cite[p. 16]{AkhiezerGlazman}. According to \cite[formula
7.414 (2)]{GR},
\[
\int_{0}^{\infty}L_{n}(t)e^{-\left(  \frac{1}{2}+i\rho\right)  t}%
dt=\frac{\left(  -1\right)  ^{n}\left(  \frac{1}{2}-i\rho\right)  ^{n}%
}{\left(  \frac{1}{2}+i\rho\right)  ^{n+1}}.
\]
Hence
\[
\widetilde{\psi}(\rho,x)=e^{-i\rho x}\left(  \left(
\begin{tabular}
[c]{c}%
$1$\\
$0$%
\end{tabular}
\ \right)  +\sum_{n=0}^{\infty}\frac{\left(  -1\right)  ^{n}\left(  \frac
{1}{2}-i\rho\right)  ^{n}}{\left(  \frac{1}{2}+i\rho\right)  ^{n+1}}%
\widetilde{a}_{n}(x)\right)  .
\]
Note that
\[
\frac{1}{2}+i\rho=\frac{1}{1+\widetilde{z}}.
\]
Thus, substitution of $\widetilde{z}$ leads to (\ref{psitil=}). Notice that
$z\left(  \rho\right)  $ is a M\"{o}bius transformation of the upper halfplane
of the complex variable $\rho$ onto the unit disk $D=\left\{  z\in
\mathbb{C}:\,\left\vert z\right\vert \leq1\right\}  $, and $\widetilde
{z}\left(  \rho\right)  $ transforms the lower halfplane onto the unit disk.
In particular, since $\widetilde{\psi}_{1,2}(\rho,x)e^{i\rho x}$ are analytic
functions in the lower halfplane $\operatorname{Im}\rho<0$, after the
M\"{o}bius transformation $\rho\longmapsto\widetilde{z}\left(  \rho\right)  $
they become analytic in the unit disk $\widetilde{D}$, and their power series
in $\widetilde{z}$ converge in $\widetilde{D}$. The proof of the series
representations for the other three Jost solutions is analogous.$\qquad
\square$

\begin{remark}
For any $x\in\mathbb{R}$ the series $\sum_{n=0}^{\infty}\left\vert
a_{1,n}\right\vert ^{2}(x)$, $\sum_{n=0}^{\infty}\left\vert a_{2,n}\right\vert
^{2}(x)$, as well as the corresponding series of the other coefficients from
(\ref{A}), (\ref{B}) converge. This is a consequence of the fact that
$a_{1,n}$, $a_{2,n}$ etc. are Fourier coefficients of functions from
$\mathcal{L}_{2}\left(  0,\infty;e^{-t}\right)  $ with respect to the Laguerre
polynomials. Hence, for any $x\in\mathbb{R}$ the functions $\psi_{1,2}%
(\rho,x)e^{-i\rho x}$, $\varphi_{1,2}(\rho,x)e^{i\rho x}$ belong to the Hardy
space of the unit disc $H^{2}(D)$ as functions of $z$ (see, e.g.,
\cite[Theorem 17.12]{Rudin}). Simlarly, the functions $\widetilde{\psi}%
_{1,2}(\rho,x)e^{i\rho x}$, $\widetilde{\varphi}_{1,2}(\rho,x)e^{-i\rho x}$
belong to the Hardy space of the unit disc $H^{2}(\widetilde{D})$ as functions
of $\widetilde{z}$.
\end{remark}

Since the series in (\ref{phi=})-(\ref{psitil=}) are power series with respect
to the parameters $z$ or $\widetilde{z}$, respectively, which are related to
the spectral parameter $\rho$, we refer to the series representations
(\ref{phi=})-(\ref{psitil=}) of the Jost solutions as to the \textbf{spectral
parameter power series} (SPPS), in analogy with the SPPS introduced in
\cite{KrPorter2010} for regular solutions of Sturm-Liouville equations.

\section{Recurrent integration procedure for the
coefficients\label{Sect recurrent}}

In this section we devise the formulas for the calculation of the SPPS
coefficients introduced in the previous section. We begin by obtaining
expressions for the first coefficients.

\subsection{First coefficients}

Since $z\left(  \frac{i}{2}\right)  =0$, from the series representations
(\ref{phi=}) and (\ref{psi=}) we obtain%
\begin{equation}
\varphi(\frac{i}{2},x)=e^{\frac{x}{2}}\left(  \left(
\begin{tabular}
[c]{c}%
$1$\\
$0$%
\end{tabular}
\ \right)  +b_{0}(x)\right)  ,\quad\psi(\frac{i}{2},x)=e^{-\frac{x}{2}}\left(
\left(
\begin{tabular}
[c]{c}%
$0$\\
$1$%
\end{tabular}
\ \right)  +a_{0}(x)\right)  . \label{phi(i/2)}%
\end{equation}
Similarly, since $\widetilde{z}\left(  -\frac{i}{2}\right)  =0$, from the
series representations (\ref{phitil=}) and (\ref{psitil=}) we obtain
\[
\widetilde{\varphi}(-\frac{i}{2},x)=e^{\frac{x}{2}}\left(  \left(
\begin{tabular}
[c]{c}%
$0$\\
$-1$%
\end{tabular}
\ \right)  +\widetilde{b}_{0}(x)\right)  ,\quad\widetilde{\psi}(-\frac{i}%
{2},x)=e^{-\frac{x}{2}}\left(  \left(
\begin{tabular}
[c]{c}%
$1$\\
$0$%
\end{tabular}
\ \right)  +\widetilde{a}_{0}(x)\right)  .
\]

Thus, to obtain the first SPPS coefficients, it is sufficient to compute the
Jost solutions $\varphi$ and $\psi$ corresponding to $\rho=\frac{i}{2}$ and
the Jost solutions $\widetilde{\varphi}$ and $\widetilde{\psi}$ corresponding
to $\rho=-\frac{i}{2}$. Then
\begin{equation}
b_{0}(x)=e^{-\frac{x}{2}}\varphi(\frac{i}{2},x)-\left(
\begin{tabular}
[c]{c}%
$1$\\
$0$%
\end{tabular}
\right)  ,\quad\widetilde{b}_{0}(x)=e^{-\frac{x}{2}}\widetilde{\varphi}%
(-\frac{i}{2},x)-\left(
\begin{tabular}
[c]{c}%
$0$\\
$-1$%
\end{tabular}
\right)  , \label{b0}%
\end{equation}%
\begin{equation}
a_{0}(x)=e^{\frac{x}{2}}\psi(\frac{i}{2},x)-\left(
\begin{tabular}
[c]{c}%
$0$\\
$1$%
\end{tabular}
\right)  ,\quad\widetilde{a}_{0}(x)=e^{\frac{x}{2}}\widetilde{\psi}(-\frac
{i}{2},x)-\left(
\begin{tabular}
[c]{c}%
$1$\\
$0$%
\end{tabular}
\right)  . \label{a0}%
\end{equation}
To obtain a procedure for computing the rest of the SPPS coefficients, let us
substitute the series representations (\ref{phi=})-(\ref{psitil=}) into
(\ref{ZS1}), (\ref{ZS2}).

\subsection{Coefficients $a_{n}(x)$}

Substitution of (\ref{psi=}) leads to the equations%
\begin{equation}
\left(  z+1\right)  \sum_{n=0}^{\infty}\left(  -1\right)  ^{n}z^{n}%
a_{1,n}^{\prime}(x)+\left(  z-1\right)  \sum_{n=0}^{\infty}\left(  -1\right)
^{n}z^{n}a_{1,n}(x)=q(x)\left(  z+1\right)  \sum_{n=0}^{\infty}\left(
-1\right)  ^{n}z^{n}a_{2,n}(x)+q(x), \label{eq for an}%
\end{equation}%
\[
\sum_{n=0}^{\infty}\left(  -1\right)  ^{n}z^{n}a_{2,n}^{\prime}(x)=r(x)\sum
_{n=0}^{\infty}\left(  -1\right)  ^{n}z^{n}a_{1,n}(x),
\]
where to obtain the first one we took into account that
\[
2i\rho=\frac{z-1}{z+1}.
\]
Now, from the second equation, by equating the coefficients of two convergent
power series, we immediately obtain the relation
\begin{equation}
a_{2,n}^{\prime}(x)=r(x)a_{1,n}(x),\quad n=0,1,\ldots, \label{a2n prime}%
\end{equation}
while the first one gives the relations%
\begin{equation}
a_{1,0}^{\prime}(x)-a_{1,0}(x)-q(x)a_{2,0}(x)=q(x) \label{a10 prime}%
\end{equation}
and
\begin{equation}
a_{1,n}^{\prime}(x)-a_{1,n}(x)-q(x)a_{2,n}(x)=a_{1,n-1}^{\prime}%
(x)+a_{1,n-1}(x)-q(x)a_{2,n-1}(x),\quad n=1,2,\ldots. \label{a1n prime}%
\end{equation}
From (\ref{a2n prime}) we have
\begin{equation}
a_{1,n}(x)=\frac{a_{2,n}^{\prime}(x)}{r(x)},\quad n=0,1,\ldots. \label{a1n=}%
\end{equation}
Substitution of this relation into (\ref{a1n prime}) leads to the equations%
\[
\left(  \frac{d}{dx}-1\right)  \left(  \frac{1}{r(x)}\frac{d}{dx}%
a_{2,n}(x)\right)  -q(x)a_{2,n}(x)=a_{1,n-1}^{\prime}(x)+a_{1,n-1}%
(x)-q(x)a_{2,n-1}(x),\quad n=1,2,\ldots.
\]
Note that $\left(  \frac{d}{dx}-1\right)  \cdot\,=e^{x}\frac{d}{dx}\left(
e^{-x}\cdot\right)  $. Thus,%
\begin{equation}
\frac{d}{dx}\left(  R(x)\frac{d}{dx}a_{2,n}(x)\right)  -Q(x)a_{2,n}%
(x)=h_{n}(x),\quad n=1,2,\ldots, \label{dRd}%
\end{equation}
where%
\begin{equation}
R(x):=\frac{e^{-x}}{r(x)},\quad Q(x):=e^{-x}q(x) \label{R and Q}%
\end{equation}
and
\[
h_{n}(x):=e^{-x}\left(  a_{1,n-1}^{\prime}(x)+a_{1,n-1}(x)-q(x)a_{2,n-1}%
(x)\right)  .
\]
Note that substitution of $a_{1,0}(x)=\frac{a_{2,0}^{\prime}(x)}{r(x)}$ into
(\ref{a10 prime}) leads to the equation
\[
\frac{d}{dx}\left(  R(x)\frac{d}{dx}a_{2,0}(x)\right)  -Q(x)a_{2,0}(x)=Q(x),
\]
which for $f(x):=a_{2,0}(x)+1=e^{\frac{x}{2}}\psi_{2}(\frac{i}{2},x)$ takes
the form%
\begin{equation}
\frac{d}{dx}\left(  R(x)\frac{d}{dx}f(x)\right)  -Q(x)f(x)=0. \label{eq for f}%
\end{equation}
Thus, $f(x)$ satisfies this homogeneous second order linear equation together
with the asymptotic relation $f(x)\sim1$, $x\rightarrow+\infty$. For
simplicity, let us assume that $f(x)$ is a nonvanishing function. If, on the
contrary, it has zeros, then with the aid of $f(x)$ one can construct another
linearly independent solution of (\ref{eq for f}) and consider a nonvanishing
linear combination of two solutions. We will not dwell on this issue, assuming
that $f(x)$ is nonvanishing.

With the aid of the function $f(x)$, equation (\ref{dRd}) can be written in
the form
\[
\frac{1}{f(x)}\frac{d}{dx}\left(  R(x)f^{2}(x)\frac{d}{dx}\left(
\frac{a_{2,n}(x)}{f(x)}\right)  \right)  =h_{n}(x),\quad n=1,2,\ldots,
\]
which leads to a simple representation for its solution
\begin{equation}
a_{2,n}(x)=f(x)\int_{x}^{\infty}\frac{1}{R(t)f^{2}(t)}\int_{t}^{\infty
}f(s)h_{n}(s)dsdt,\quad n=1,2,\ldots. \label{a1}%
\end{equation}

Note that for calculating $a_{2,n}(x)$ as well as $a_{1,n}(x)$, $a_{1,n}%
^{\prime}(x)$, $a_{2,n}^{\prime}(x)$ no differentiation is required apart from
$f^{\prime}(x)$, which needs to be computed once. Indeed, we have
\begin{equation}
a_{2,n}^{\prime}(x)=\frac{f^{\prime}(x)}{f(x)}a_{2,n}(x)-\frac{e^{x}%
r(x)}{f(x)}\int_{x}^{\infty}f(s)h_{n}(s)ds,\quad n=1,2,\ldots,\label{a2}%
\end{equation}%
\begin{equation}
a_{1,n}(x)=\frac{a_{2,n}^{\prime}(x)}{r(x)}=\frac{f^{\prime}(x)}%
{r(x)f(x)}a_{2,n}(x)-\frac{e^{x}}{f(x)}\int_{x}^{\infty}f(s)h_{n}(s)ds,\quad
n=1,2,\ldots,\label{a3}%
\end{equation}
and from (\ref{a1n prime}):
\begin{equation}
a_{1,n}^{\prime}(x)=a_{1,n-1}^{\prime}(x)+a_{1,n}(x)+a_{1,n-1}(x)+q(x)(a_{2,n}%
(x)-a_{2,n-1}(x)),\quad n=1,2,\ldots,\label{a4}%
\end{equation}
with $a_{1,0}^{\prime}(x)=a_{1,0}(x)+q(x)(a_{2,0}(x)+1)$ (see (\ref{a10 prime})).

\subsection{Coefficients $\widetilde{a}_{n}(x)$}

Consider the series representation (\ref{psitil=}) for $\widetilde{\psi}%
(\rho,x)$. Its substitution into (\ref{ZS1}), (\ref{ZS2}) leads to the
equations%
\[
\widetilde{a}_{1,n}^{\prime}(x)=q(x)\widetilde{a}_{2,n}(x)
\]
and%
\[
\left(  \widetilde{z}+1\right)  \sum_{n=0}^{\infty}\left(  -1\right)
^{n}\widetilde{z}^{n}\widetilde{a}_{2,n}^{\prime}(x)+\left(  \widetilde
{z}-1\right)  \sum_{n=0}^{\infty}\left(  -1\right)  ^{n}\widetilde{z}%
^{n}\widetilde{a}_{2,n}(x)=r(x)\left(  \widetilde{z}+1\right)  \sum
_{n=0}^{\infty}\left(  -1\right)  ^{n}\widetilde{z}^{n}\widetilde{a}%
_{1,n}(x)+r(x).
\]
Their comparison with (\ref{a1n=}) and (\ref{eq for an}) shows that they
coincide up to the following change of the notations:%
\[
a_{1,n}(x)\leftrightsquigarrow\widetilde{a}_{2,n}(x),\quad a_{2,n}%
(x)\leftrightsquigarrow\widetilde{a}_{1,n}(x),\quad q(x)\leftrightsquigarrow
r(x),\quad r(x)\leftrightsquigarrow q(x).
\]
This gives us immediately the recurrent integration formulas for the
coefficients $\widetilde{a}_{1,n}(x)$, $\widetilde{a}_{2,n}(x)$:
\begin{equation}
\widetilde{a}_{1,n}(x)=\widetilde{f}(x)\int_{x}^{\infty}\frac{e^{t}%
q(t)}{\widetilde{f}^{2}(t)}\int_{t}^{\infty}\widetilde{f}(s)\widetilde{h}%
_{n}(s)dsdt,\quad n=1,2,\ldots. \label{atil1}%
\end{equation}

\begin{equation}
\widetilde{a}_{1,n}^{\prime}(x)=\frac{\widetilde{f}^{\prime}(x)}{\widetilde
{f}(x)}\widetilde{a}_{1,n}(x)-\frac{e^{x}q(x)}{\widetilde{f}(x)}\int
_{x}^{\infty}\widetilde{f}(s)\widetilde{h}_{n}(s)ds,\quad n=1,2,\ldots,
\label{atil2}%
\end{equation}%
\begin{equation}
\widetilde{a}_{2,n}(x)=\frac{\widetilde{a}_{1,n}^{\prime}(x)}{q(x)}%
=\frac{\widetilde{f}^{\prime}(x)}{q(x)\widetilde{f}(x)}\widetilde{a}%
_{1,n}(x)-\frac{e^{x}}{\widetilde{f}(x)}\int_{x}^{\infty}\widetilde
{f}(s)\widetilde{h}_{n}(s)ds,\quad n=1,2,\ldots, \label{atil3}%
\end{equation}%
\begin{equation}
\widetilde{a}_{2,n}^{\prime}(x)=\widetilde{a}_{2,n-1}^{\prime}(x)+\widetilde
{a}_{2,n}(x)+\widetilde{a}_{2,n-1}(x)+r(x)(\widetilde{a}_{1,n}(x)-\widetilde
{a}_{1,n-1}(x)),\quad n=1,2,\ldots, \label{atil4}%
\end{equation}
with $\widetilde{a}_{2,0}^{\prime}(x)=\widetilde{a}_{2,0}(x)+r(x)(\widetilde
{a}_{1,0}(x)+1)$. Here%
\[
\widetilde{h}_{n}(x)=e^{-x}\left(  \widetilde{a}_{2,n-1}^{\prime
}(x)+\widetilde{a}_{2,n-1}(x)-r(x)\widetilde{a}_{1,n-1}(x)\right)  ,
\]
and $\widetilde{f}$ is a solution of
\[
\frac{d}{dx}\left(  \widetilde{Q}(x)\frac{d}{dx}\widetilde{f}(x)\right)
-\widetilde{R}(x)\widetilde{f}(x)=0,\quad\widetilde{f}(x)\sim1,\,x\rightarrow
+\infty,
\]
$\widetilde{f}(x):=\widetilde{a}_{1,0}(x)+1=e^{\frac{x}{2}}\widetilde{\psi
}_{1}(-\frac{i}{2},x)$, $\widetilde{Q}(x):=\frac{e^{-x}}{q(x)}$,
$\widetilde{R}(x):=e^{-x}r(x)$.

\subsection{Coefficients $b_{n}(x)$}

Consider the series representation (\ref{phi=}) for $\varphi(\rho,x)$. Its
substitution into (\ref{ZS1}), (\ref{ZS2}) leads to the equations%
\begin{equation}
b_{1,n}^{\prime}(x)=q(x)b_{2,n}(x),\quad n=0,1,\ldots, \label{b1nPrime}%
\end{equation}%
\[
b_{2,0}^{\prime}(x)+b_{2,0}(x)=r(x)(b_{1,0}(x)+1),
\]%
\[
b_{2,n}^{\prime}(x)+b_{2,n}(x)-r(x)b_{1,n}(x)=b_{2,n-1}^{\prime}%
(x)-b_{2,n-1}(x)-r(x)b_{1,n-1}(x),\quad n=1,2,\ldots.
\]
Substituting $b_{2,n}(x)=b_{1,n}^{\prime}(x)/q(x)$ (from (\ref{b1nPrime}))
into the last equation, we obtain (compare with (\ref{dRd}))%
\[
\frac{d}{dx}\left(  \frac{1}{Q(x)}\frac{d}{dx}b_{1,n}(x)\right)  -\frac
{1}{R(x)}b_{1,n}(x)=p_{n}(x),\quad n=1,2,\ldots,
\]
where $Q(x)$ and $R(x)$ are those from (\ref{R and Q}), and
\[
p_{n}(x):=e^{x}\left(  b_{2,n-1}^{\prime}(x)-b_{2,n-1}(x)-r(x)b_{1,n-1}%
(x)\right)  .
\]
For $n=0$ we have the equation for $g(x):=b_{1,0}(x)+1=e^{-\frac{x}{2}}%
\varphi_{1}(\frac{i}{2},x)$:
\[
\frac{d}{dx}\left(  \frac{1}{Q(x)}\frac{d}{dx}g(x)\right)  -\frac{1}%
{R(x)}g(x)=0
\]
and $g(x)\sim1$, $x\rightarrow-\infty$ (compare with (\ref{eq for f})).

By analogy with the previous cases we obtain the formulas of recurrent
integration procedure:%
\begin{equation}
b_{1,n}(x)=g(x)\int_{-\infty}^{x}\frac{Q(t)}{g^{2}(t)}\int_{-\infty}%
^{t}g(s)p_{n}(s)dsdt, \label{b1}%
\end{equation}%
\begin{equation}
b_{1,n}^{\prime}(x)=\frac{g^{\prime}(x)}{g(x)}b_{1,n}(x)+\frac{Q(x)}{g(x)}%
\int_{-\infty}^{x}g(s)p_{n}(s)ds, \label{b2}%
\end{equation}%
\begin{equation}
b_{2,n}(x)=\frac{b_{1,n}^{\prime}(x)}{q(x)}=\frac{g^{\prime}(x)}%
{q(x)g(x)}b_{1,n}(x)+\frac{e^{-x}}{g(x)}\int_{-\infty}^{x}g(s)p_{n}(s)ds,
\label{b3}%
\end{equation}%
\begin{equation}
b_{2,n}^{\prime}(x)=b_{2,n-1}^{\prime}(x)-b_{2,n}(x)-b_{2,n-1}(x)+r(x)\left(
b_{1,n}(x)-b_{1,n-1}(x)\right)  \label{b4}%
\end{equation}
with%
\[
b_{2,0}^{\prime}(x)=-b_{2,0}(x)+r(x)\left(  b_{1,0}(x)+1\right)  .
\]

\subsection{Coefficients $\widetilde{b}_{n}(x)$}

Finally, consider $\widetilde{\varphi}(\rho,x)$. Substitution of
(\ref{phitil=}) into (\ref{ZS1}), (\ref{ZS2}) leads to the equations%
\begin{equation}
\widetilde{b}_{2,n}^{\prime}(x)=r(x)\widetilde{b}_{1,n}(x),\quad
n=0,1,\ldots,\label{btil2n}%
\end{equation}%
\begin{equation}
\widetilde{b}_{1,0}^{\prime}(x)+\widetilde{b}_{1,0}(x)=q(x)\left(
\widetilde{b}_{2,0}(x)-1\right)  ,\label{btil10}%
\end{equation}%
\begin{equation}
\widetilde{b}_{1,n}^{\prime}(x)+\widetilde{b}_{1,n}(x)-q(x)\widetilde{b}%
_{2,n}(x)=\widetilde{b}_{1,n-1}^{\prime}(x)-\widetilde{b}_{1,n-1}%
(x)-q(x)\widetilde{b}_{2,n-1}(x).\label{b1tiln}%
\end{equation}
Substitution of $\widetilde{b}_{1,n}(x)=\widetilde{b}_{2,n}^{\prime}(x)/r(x)$
(from (\ref{btil2n})) into (\ref{btil10}) and (\ref{b1tiln}) leads to the
equations%
\[
\frac{d}{dx}\left(  \frac{1}{\widetilde{R}(x)}\frac{d}{dx}\widetilde
{g}(x)\right)  -\frac{1}{\widetilde{Q}(x)}\widetilde{g}(x)=0,
\]%
\[
\frac{d}{dx}\left(  \frac{1}{\widetilde{R}(x)}\frac{d}{dx}\widetilde{b}%
_{2,n}(x)\right)  -\frac{1}{\widetilde{Q}(x)}\widetilde{b}_{2,n}%
(x)=\widetilde{p}_{n}(x),
\]
where $\widetilde{g}(x):=\widetilde{b}_{2,0}(x)-1=e^{-\frac{x}{2}}%
\widetilde{\varphi}_{2}(-\frac{i}{2},x)$ (note that $\widetilde{g}(x)\sim-1$,
$x\rightarrow-\infty$),
\[
\widetilde{p}_{n}(x):=e^{x}\left(  \widetilde{b}_{1,n-1}^{\prime
}(x)-\widetilde{b}_{1,n-1}(x)-q(x)\widetilde{b}_{2,n-1}(x)\right)
\]
(compare with (\ref{eq for f}) and (\ref{dRd})). Similarly to the previous
cases we obtain the recurrent integration formulas for the coefficients and
their derivatives.
\begin{equation}
\widetilde{b}_{2,n}(x)=\widetilde{g}(x)\int_{-\infty}^{x}\frac{\widetilde
{R}(t)}{\widetilde{g}^{2}(t)}\int_{-\infty}^{t}\widetilde{g}(s)\widetilde
{p}_{n}(s)dsdt,\label{btil1}%
\end{equation}%
\begin{equation}
\widetilde{b}_{2,n}^{\prime}(x)=\frac{\widetilde{g}^{\prime}(x)}{\widetilde
{g}(x)}\widetilde{b}_{2,n}(x)+\frac{\widetilde{R}(x)}{\widetilde{g}(x)}%
\int_{-\infty}^{x}\widetilde{g}(s)\widetilde{p}_{n}(s)ds,\label{btil2}%
\end{equation}%
\begin{equation}
\widetilde{b}_{1,n}(x)=\frac{\widetilde{b}_{2,n}^{\prime}(x)}{r(x)}%
=\frac{\widetilde{g}^{\prime}(x)}{r(x)\widetilde{g}(x)}\widetilde{b}%
_{2,n}(x)+\frac{e^{-x}}{\widetilde{g}(x)}\int_{-\infty}^{x}\widetilde
{g}(s)\widetilde{p}_{n}(s)ds,\label{btil3}%
\end{equation}%
\begin{equation}
\widetilde{b}_{1,n}^{\prime}(x)=\widetilde{b}_{1,n-1}^{\prime}(x)-\widetilde
{b}_{1,n}(x)-\widetilde{b}_{1,n-1}(x)+q(x)\left(  \widetilde{b}_{2,n}%
(x)-\widetilde{b}_{2,n-1}(x)\right)  \label{btil4}%
\end{equation}
with
\[
\widetilde{b}_{1,0}^{\prime}(x)=-\widetilde{b}_{1,0}(x)+q(x)\left(
\widetilde{b}_{2,0}(x)-1\right)  .
\]

\subsection{Asymptotics of the coefficients}

The next statement justifies the integration in the recurrent procedure for
constructing the coefficients. We prove it for the coefficients of $\psi
(\rho,x)$. In the other three cases the proof is analogous.

\begin{proposition}
Let $q$, $r$ $\in\mathcal{P}$. Then the coefficients $a_{1,n}(x)$,
$a_{2,n}(x)$ and their derivatives satisfy the asymptotic relations for
$x\rightarrow+\infty$:%
\begin{equation}
a_{1,n}(x)=o(1),\quad a_{1,n}^{\prime}(x)=o(1),\quad n=0,1,\ldots,
\label{asymp a1}%
\end{equation}%
\begin{equation}
a_{2,0}(x)=o(1),\quad a_{2,n}(x)=o(\int_{x}^{\infty}r(x)dx),\quad
n=1,2,\ldots, \label{asympt a2}%
\end{equation}%
\begin{equation}
a_{2,n}^{\prime}(x)=o(r(x)),\quad n=0,1,\ldots. \label{asympt a2prime}%
\end{equation}

\end{proposition}

\textbf{Proof. }We prove first the relations for $n=0$. Consider
\[
a_{2,0}(x)=f(x)-1=e^{\frac{x}{2}}\psi_{2}(\frac{i}{2},x)-1=e^{\frac{x}{2}%
}e^{-\frac{x}{2}}\left(  1+o(1)\right)  -1=o(1),\quad x\rightarrow+\infty.
\]
Next,
\begin{align*}
a_{2,0}^{\prime}(x)  &  =f^{\prime}(x)=\frac{f(x)}{2}+e^{\frac{x}{2}}\psi
_{2}^{\prime}(\frac{i}{2},x)\\
&  =\frac{f(x)}{2}+e^{\frac{x}{2}}\left(  r(x)\psi_{1}(\frac{i}{2}%
,x)-\frac{\psi_{2}(\frac{i}{2},x)}{2}\right)  ,
\end{align*}
where we used the fact that $\psi_{1}(\frac{i}{2},x)$, $\psi_{2}(\frac{i}%
{2},x)$ are solutions of (\ref{ZS1}), (\ref{ZS2}) with $\rho=\frac{i}{2}$.
Thus,%
\[
a_{2,0}^{\prime}(x)=e^{\frac{x}{2}}r(x)\psi_{1}(\frac{i}{2},x)=r(x)e^{\frac
{x}{2}}e^{-\frac{x}{2}}o(1)=o(r(x)),\quad x\rightarrow+\infty.
\]
Now, from (\ref{a1n=}) we obtain
\begin{equation}
a_{1,0}(x)=\frac{a_{2,0}^{\prime}(x)}{r(x)}=\frac{r(x)o(1)}{r(x)}=o(1),\quad
x\rightarrow+\infty, \label{a10prime}%
\end{equation}
and from (\ref{a10 prime}):%
\[
a_{1,0}^{\prime}(x)=a_{1,0}(x)+q(x)(a_{2,0}(x)+1)=o(1)+q(x)f(x)=o(1),\quad
x\rightarrow+\infty,
\]
because $q(x)=o(1)$ and $f(x)\sim1$. Thus, the asymptotic relations for $n=0$
are proved. For $n=1,2,\ldots$ we prove them by induction. First, let us prove
the asymptotic relations for $n=1$. Consider%
\begin{equation}
a_{2,1}(x)=f(x)\int_{x}^{\infty}\frac{e^{t}r(t)}{f^{2}(t)}\int_{t}^{\infty
}f(s)h_{1}(s)dsdt. \label{a21}%
\end{equation}
Here
\[
h_{1}(x)=e^{-x}\left(  a_{1,0}^{\prime}(x)+a_{1,0}(x)-q(x)a_{2,0}(x)\right)
=e^{-x}o(1),\quad x\rightarrow+\infty
\]
(we used the previously proved asymptotic). Hence $\int_{t}^{\infty}%
f(s)h_{1}(s)ds=o(e^{-t})$, and substitution of this asymptotic into
(\ref{a21}) leads to
\[
a_{2,1}(x)=o(\int_{x}^{\infty}r(x)dx).
\]
Next,
\begin{align}
a_{2,1}^{\prime}(x)  &  =\frac{f^{\prime}(x)}{f(x)}a_{2,1}(x)-\frac{e^{x}%
r(x)}{f(x)}\int_{x}^{\infty}f(t)h_{1}(t)dt\nonumber\\
&  =o(r(x))o(\int_{x}^{\infty}r(x)dx)-e^{x}r(x)o(e^{-x})\nonumber\\
&  =o(r(x)). \label{a21prime}%
\end{align}
The asymptotic for $a_{1,1}(x)$ and $a_{1,1}^{\prime}(x)$ is proved
analogously to that of $a_{1,0}(x)$ and $a_{1,0}^{\prime}(x)$ above. Now,
assuming the validity of the asymptotic relations for $n-1$, we prove them for
$n$. Consider%
\[
a_{2,n}(x)=f(x)\int_{x}^{\infty}\frac{e^{t}r(t)}{f^{2}(t)}\int_{t}^{\infty
}f(s)h_{n}(s)dsdt.
\]
We have $h_{n}(x)=e^{-x}\left(  a_{1,n-1}^{\prime}(x)+a_{1,n-1}%
(x)-q(x)a_{2,n-1}(x)\right)  =e^{-x}o(1)$, $x\rightarrow+\infty$, and hence
substituting this into the expression for $a_{2,n}(x)$ we obtain that
(\ref{asympt a2}) is valid. Next, the asymptotic for $a_{2,n}^{\prime}(x)$ is
proved analogously to that of $a_{2,1}^{\prime}(x)$, see (\ref{a21prime}), as
well as the asymptotic of $a_{1,n}(x)$ analogously to that of $a_{1,0}(x)$,
see (\ref{a10prime}). Finally,
\begin{align*}
a_{1,n}^{\prime}(x)  &  =a_{1,n-1}^{\prime}(x)+a_{1,n}(x)+a_{1,n-1}%
(x)+q(x)(a_{2,n}(x)-a_{2,n-1}(x))\\
&  =o(1)+o(1)+o(1)+o(1)o(\int_{x}^{\infty}r(x)dx)=o(1),\quad x\rightarrow
+\infty,
\end{align*}
that finishes the proof.$\qquad\square$

\subsection{Summary of the recurrent
procedure\label{subsect summary recurrent}}

Let us summarize the procedure for computing the coefficients of the series
(\ref{phi=})-(\ref{psitil=}). First, the Jost solutions $\varphi(\frac{i}%
{2},x)$, $\psi(\frac{i}{2},x)$, $\widetilde{\varphi}(-\frac{i}{2},x)$,
$\widetilde{\psi}(-\frac{i}{2},x)$ need to be computed, that gives us
\[
f(x)=a_{2,0}(x)+1=e^{\frac{x}{2}}\psi_{2}(\frac{i}{2},x),\quad f^{\prime
}(x)=a_{2,0}^{\prime}(x)=e^{\frac{x}{2}}r(x)\psi_{1}(\frac{i}{2},x),
\]
\[
\widetilde{f}(x)=\widetilde{a}_{1,0}(x)+1=e^{\frac{x}{2}}\widetilde{\psi}%
_{1}(-\frac{i}{2},x),\quad\widetilde{f}^{\prime}(x)=\widetilde{a}%
_{1,0}^{\prime}(x)=e^{\frac{x}{2}}q(x)\widetilde{\psi}_{2}(-\frac{i}{2},x),
\]%
\[
g(x)=b_{1,0}(x)+1=e^{-\frac{x}{2}}\varphi_{1}(\frac{i}{2},x),\quad g^{\prime
}(x)=b_{1,0}^{\prime}(x)=e^{-\frac{x}{2}}q(x)\varphi_{2}(\frac{i}{2},x),
\]%
\[
\widetilde{g}(x)=\widetilde{b}_{2,0}(x)-1=e^{-\frac{x}{2}}\widetilde{\varphi
}_{2}(-\frac{i}{2},x),\quad\widetilde{g}^{\prime}(x)=\widetilde{b}%
_{2,0}^{\prime}(x)=e^{-\frac{x}{2}}r(x)\widetilde{\varphi}_{1}(-\frac{i}%
{2},x).
\]
The rest of the coefficients are calculated with the aid of the formulas
(\ref{a1})-(\ref{a4}) for $\left\{  a_{n}(x)\right\}  $, (\ref{atil1}%
)-(\ref{atil4}) for $\left\{  \widetilde{a}_{n}(x)\right\}  $, (\ref{b1}%
)-(\ref{b4}) for $\left\{  b_{n}(x)\right\}  $ and (\ref{btil1})-(\ref{btil4})
for $\left\{  \widetilde{b}_{n}(x)\right\}  $.

Note that $f$ and $f^{\prime}$ can be computed by solving (\ref{eq for f})
with the boundary conditions $f(x_{0})=1$ and $f^{\prime}(x_{0})=0$ at a
sufficiently large $x_{0}$. This can be done with the aid of any numerical
method, e.g., the SPPS method from \cite{KrPorter2010}, as it is done in the
present work. The other three functions $\widetilde{f}$, $g$, $\widetilde{g}$
and their derivatives can be computed in a similar way.

\section{Solution of the direct problem\label{SectDirect}}

Substitution of the SPPS (\ref{phi=})-(\ref{psitil=}) into (\ref{atila}) leads
to the SPPS representations for $\mathbf{a}(\rho)$ and $\widetilde{\mathbf{a}%
}(\rho)$:%
\begin{align}
\mathbf{a}(\rho)  &  =\varphi_{1}(\rho,0)\psi_{2}(\rho,0)-\varphi_{2}%
(\rho,0)\psi_{1}(\rho,0)\nonumber\\
&  =\left(  1+\left(  z+1\right)  \sum_{n=0}^{\infty}\left(  -1\right)
^{n}z^{n}b_{1,n}(0)\right)  \left(  1+\left(  z+1\right)  \sum_{n=0}^{\infty
}\left(  -1\right)  ^{n}z^{n}a_{2,n}(0)\right) \nonumber\\
&  -\left(  z+1\right)  ^{2}\left(  \sum_{n=0}^{\infty}\left(  -1\right)
^{n}z^{n}b_{2,n}(0)\right)  \left(  \sum_{n=0}^{\infty}\left(  -1\right)
^{n}z^{n}a_{1,n}(0)\right)  \label{SPPS a}%
\end{align}
for all $\rho\in\overline{\mathbb{C}^{+}}$ and
\begin{align}
\widetilde{\mathbf{a}}(\rho)  &  =\widetilde{\varphi}_{1}(\rho,0)\widetilde
{\psi}_{2}(\rho,0)-\widetilde{\varphi}_{2}(\rho,0)\widetilde{\psi}_{1}%
(\rho,0)\nonumber\\
&  =\left(  \widetilde{z}+1\right)  ^{2}\left(  \sum_{n=0}^{\infty}\left(
-1\right)  ^{n}\widetilde{z}^{n}\widetilde{b}_{1,n}(0)\right)  \left(
\sum_{n=0}^{\infty}\left(  -1\right)  ^{n}\widetilde{z}^{n}\widetilde{a}%
_{2,n}(0)\right) \nonumber\\
&  +\left(  1-\left(  \widetilde{z}+1\right)  \sum_{n=0}^{\infty}\left(
-1\right)  ^{n}\widetilde{z}^{n}\widetilde{b}_{2,n}(0)\right)  \left(
1+\left(  \widetilde{z}+1\right)  \sum_{n=0}^{\infty}\left(  -1\right)
^{n}\widetilde{z}^{n}\widetilde{a}_{1,n}(0)\right)  \label{SPPS atil}%
\end{align}
for all $\rho\in\overline{\mathbb{C}^{-}}$.

Analogously, the SPPS representations for $\mathbf{b}(\rho)$ and
$\widetilde{\mathbf{b}}(\rho)$ are obtained by substituting (\ref{phi=}%
)-(\ref{psitil=}) into (\ref{btilb}):%
\begin{align}
\mathbf{b}(\rho)  &  =\varphi_{2}(\rho,0)\widetilde{\psi}_{1}(\rho
,0)-\varphi_{1}(\rho,0)\widetilde{\psi}_{2}(\rho,0)\nonumber\\
&  =\left(  z+1\right)  \left(  \sum_{n=0}^{\infty}\left(  -1\right)
^{n}z^{n}b_{2,n}(0)\right)  \left(  1+\left(  \widetilde{z}+1\right)
\sum_{n=0}^{\infty}\left(  -1\right)  ^{n}\widetilde{z}^{n}\widetilde{a}%
_{1,n}(0)\right) \nonumber\\
&  -\left(  \widetilde{z}+1\right)  \left(  1+\left(  z+1\right)  \sum
_{n=0}^{\infty}\left(  -1\right)  ^{n}z^{n}b_{1,n}(0)\right)  \left(
\sum_{n=0}^{\infty}\left(  -1\right)  ^{n}\widetilde{z}^{n}\widetilde{a}%
_{2,n}(0)\right)  \label{SPPS b}%
\end{align}
and
\begin{align}
\widetilde{\mathbf{b}}(\rho)  &  =\widetilde{\varphi}_{1}(\rho,0)\psi_{2}%
(\rho,0)-\widetilde{\varphi}_{2}(\rho,0)\psi_{1}(\rho,0)\nonumber\\
&  =\left(  \widetilde{z}+1\right)  \left(  \sum_{n=0}^{\infty}\left(
-1\right)  ^{n}\widetilde{z}^{n}\widetilde{b}_{1,n}(0)\right)  \left(
1+\left(  z+1\right)  \sum_{n=0}^{\infty}\left(  -1\right)  ^{n}z^{n}%
a_{2,n}(0)\right) \nonumber\\
&  +\left(  z+1\right)  \left(  1-\left(  \widetilde{z}+1\right)  \sum
_{n=0}^{\infty}\left(  -1\right)  ^{n}\widetilde{z}^{n}\widetilde{b}%
_{2,n}(0)\right)  \left(  \sum_{n=0}^{\infty}\left(  -1\right)  ^{n}%
z^{n}a_{1,n}(0)\right)  \label{SPPS btil}%
\end{align}
for all $\rho\in\mathbb{R}$. Truncating the series in (\ref{SPPS a}%
)-(\ref{SPPS btil}) up to a sufficiently large number $N$ gives us a simple
way of computing the scattering matrix elements. Note that when $\rho
\in\mathbb{R}$ both $z$ and $\widetilde{z}$ belong to the unit circle centered
at the origin; when $\rho\in\mathbb{C}^{+}$, $z$ belongs to the unit disk $D$;
when $\rho\in\mathbb{C}^{-}$, $\widetilde{z}$ belongs to the unit disk
$\widetilde{D}$. This observation implies that the computation of the
eigenvalues in $\mathbb{C}^{+}$ reduces to the location of zeros of
(\ref{SPPS a}) inside $D$, while the computation of the eigenvalues in
$\mathbb{C}^{-}$ reduces to the location of zeros of (\ref{SPPS atil}) inside
$\widetilde{D}$. Numerically, the truncated series (\ref{SPPS a}) and
(\ref{SPPS atil}) are considered, which are nothing but polynomials in the
variables $z$ and $\widetilde{z}$, respectively. Thus, the computation of the
eigenvalues reduces to the location of the roots of the respective polynomials
inside the unit disk.

The norming constants $c_{m}$ and $\widetilde{c}_{m}$ are also computed by
definition, for example as%
\[
c_{m}=\frac{\varphi_{1}(\rho_{m},0)}{\psi_{1}(\rho_{m},0)}=\frac{1+\left(
z_{m}+1\right)  \sum_{n=0}^{\infty}\left(  -1\right)  ^{n}z_{m}^{n}b_{1,n}%
(0)}{\left(  z_{m}+1\right)  \sum_{n=0}^{\infty}\left(  -1\right)  ^{n}%
z_{m}^{n}a_{1,n}(0)},
\]%
\[
\widetilde{c}_{m}=\frac{\widetilde{\varphi}_{1}(\widetilde{\rho}_{m}%
,0)}{\widetilde{\psi}_{1}(\widetilde{\rho}_{m},0)}=\frac{\left(  \widetilde
{z}_{m}+1\right)  \left(  \sum_{n=0}^{\infty}\left(  -1\right)  ^{n}%
\widetilde{z}_{m}^{n}\widetilde{b}_{1,n}(0)\right)  }{1+\left(  \widetilde
{z}_{m}+1\right)  \sum_{n=0}^{\infty}\left(  -1\right)  ^{n}\widetilde{z}%
_{m}^{n}\widetilde{a}_{1,n}(0)}.
\]
Equally, the quotients of the second components can be used. For the numerical
implementation, here again, the truncated sums should be considered.

Summarizing, the direct scattering problem is solved as follows. Given $q(x)$
and $r(x)$, compute a sufficiently large set of the coefficients $\left\{
a_{n}(0),b_{n}(0),\widetilde{a}_{n}(0),\widetilde{b}_{n}(0)\right\}
_{n=0}^{N}$ following the recurrent integration procedure from Section
\ref{Sect recurrent}, see subsection \ref{subsect summary recurrent}. Next,
the scattering matrix entries are computed for $\rho\in\mathbb{R}$ by
(\ref{SPPS a})-(\ref{SPPS btil}), where the sums truncated up to $N$ are
considered. Next, the eigenvalues $\rho_{m}\in\mathbb{C}^{+}$ are computed as
$\rho_{m}=\frac{z_{m}-1}{2i\left(  z_{m}+1\right)  }$, where $z_{m}$ are roots
of the polynomial
\begin{align*}
\mathbf{a}_{N}(z)  &  =\left(  1+\left(  z+1\right)  \sum_{n=0}^{N}\left(
-1\right)  ^{n}z^{n}b_{1,n}(0)\right)  \left(  1+\left(  z+1\right)
\sum_{n=0}^{N}\left(  -1\right)  ^{n}z^{n}a_{2,n}(0)\right) \\
&  -\left(  z+1\right)  ^{2}\left(  \sum_{n=0}^{N}\left(  -1\right)  ^{n}%
z^{n}b_{2,n}(0)\right)  \left(  \sum_{n=0}^{N}\left(  -1\right)  ^{n}%
z^{n}a_{1,n}(0)\right)
\end{align*}
located in $D$. The eigenvalues $\widetilde{\rho}_{m}\in\mathbb{C}^{-}$ are
computed as $\widetilde{\rho}_{m}=\frac{1-\widetilde{z}_{m}}{2i\left(
\widetilde{z}_{m}+1\right)  }$, where $\widetilde{z}_{m}$ are roots of the
polynomial%
\begin{align*}
\widetilde{\mathbf{a}}_{N}(\widetilde{z})  &  =\left(  \widetilde{z}+1\right)
^{2}\left(  \sum_{n=0}^{N}\left(  -1\right)  ^{n}\widetilde{z}^{n}%
\widetilde{b}_{1,n}(0)\right)  \left(  \sum_{n=0}^{N}\left(  -1\right)
^{n}\widetilde{z}^{n}\widetilde{a}_{2,n}(0)\right) \\
&  +\left(  1-\left(  \widetilde{z}+1\right)  \sum_{n=0}^{N}\left(  -1\right)
^{n}\widetilde{z}^{n}\widetilde{b}_{2,n}(0)\right)  \left(  1+\left(
\widetilde{z}+1\right)  \sum_{n=0}^{N}\left(  -1\right)  ^{n}\widetilde{z}%
^{n}\widetilde{a}_{1,n}(0)\right)
\end{align*}
located in $\widetilde{D}$.

\section{Solution of the inverse problem\label{Sect Inverse}}

First, let us notice that for recovering $q(x)$ and $r(x)$ it is sufficient to
find any of the first SPPS coefficients. For example, assume that $b_{0}(x)$
is known. This means that the Jost solution $\varphi(\frac{i}{2},x)$
corresponding to $\rho=\frac{i}{2}$ is known (formula (\ref{phi(i/2)})):
\[
\varphi_{1}(\frac{i}{2},x)=e^{\frac{x}{2}}\left(  1+b_{1,0}(x)\right)
,\quad\varphi_{2}(\frac{i}{2},x)=e^{\frac{x}{2}}b_{2,0}(x).
\]
From the AKNS system (\ref{ZS1}), (\ref{ZS2}) for $\rho=\frac{i}{2}$ we have
that%
\begin{equation}
q(x)=\frac{\varphi_{1}^{\prime}(\frac{i}{2},x)-\varphi_{1}(\frac{i}{2}%
,x)/2}{\varphi_{2}(\frac{i}{2},x)},\quad r(x)=\frac{\varphi_{2}^{\prime}%
(\frac{i}{2},x)+\varphi_{2}(\frac{i}{2},x)/2}{\varphi_{1}(\frac{i}{2}%
,x)}.\label{q and r}%
\end{equation}
Thus, knowledge of $b_{0}(x)$ is sufficient for recovering $q(x)$ and $r(x)$.
Analogously, $q(x)$ and $r(x)$ can be recovered from any of the first SPPS
coefficients. Obviously, they can be recovered as well from any linear
combination of $\varphi(\frac{i}{2},x)$ and $\psi(\frac{i}{2},x)$ or
$\widetilde{\varphi}(-\frac{i}{2},x)$ and $\widetilde{\psi}(-\frac{i}{2},x)$
which may be more convenient to avoid zeros in the denominators in
(\ref{q and r}). Let us notice that due to the asymptotic relations satisfied
by the Jost solutions, the first equality in (\ref{q and r}) may be less
convenient for recovering $q(x)$, because the absolute value of $\varphi
_{2}(\frac{i}{2},x)$ may become very small in the interval of interest (when
$x$ is negative and large), so that the division by it should be avoided.
Thus, for recovering $q(x)$, we use the equality%
\begin{equation}
q(x)=\frac{\widetilde{\varphi}_{1}^{\prime}(-\frac{i}{2},x)+\widetilde
{\varphi}_{1}(-\frac{i}{2},x)/2}{\widetilde{\varphi}_{2}(-\frac{i}{2}%
,x)},\label{q from phitil}%
\end{equation}
where $\widetilde{\varphi}_{1,2}(-\frac{i}{2},x)$ are obtained from
$\widetilde{b}_{0}(x)$. If some $x^{\ast}$ is zero of $\widetilde{\varphi}%
_{2}(-\frac{i}{2},x)$, then to compute $q(x^{\ast})$, we can use the analogous
equality
\[
q(x)=\frac{\widetilde{\psi}_{1}^{\prime}(-\frac{i}{2},x)+\widetilde{\psi}%
_{1}(-\frac{i}{2},x)/2}{\widetilde{\psi}_{2}(-\frac{i}{2},x)},
\]
where $\widetilde{\psi}_{2}(-\frac{i}{2},x^{\ast})\neq0$ since if
$\widetilde{\rho}=-\frac{i}{2}$ is not an eigenvalue $\widetilde{\varphi
}(-\frac{i}{2},x)$ and $\widetilde{\psi}(-\frac{i}{2},x)$ are linearly
independent. If  $\widetilde{\rho}=-\frac{i}{2}$ is an eigenvalue, we can use
analogous equalities for $\varphi(\frac{i}{2},x)$ and $\psi(\frac{i}{2},x)$.
Finally, if both $\rho=\frac{i}{2}$ and $\widetilde{\rho}=-\frac{i}{2}$ are
eigenvalues, and  $x^{\ast}$ is zero of both $\varphi_{2}(\frac{i}{2},x)$ and
$\widetilde{\varphi}_{2}(-\frac{i}{2},x)$, then $q(x^{\ast})$ can be computed
from (\ref{q from phitil}) as a limit value of the right-hand side expression
when $x\rightarrow x^{\ast}$.

Given the set of the scattering data (\ref{SD}), let us consider equalities
(\ref{rel 1}) and (\ref{rel 2}) at a sufficiently large number $K$ of points
$\rho_{k}\in\mathbb{R}$ and additionally, if the discrete spectrum is not
empty, the equalities (\ref{eig1}) and (\ref{eig2}) for all $\rho_{m}$,
$m=1,\ldots,M$ and $\widetilde{\rho}_{m}$, $m=1,\ldots,\widetilde{M}$. We
construct two systems of linear algebraic equations for two sets of the SPPS
coefficients $\left\{  a_{1,n}(x),b_{1,n}(x),\widetilde{a}_{1,n}%
(x),\widetilde{b}_{1,n}(x)\right\}  _{n=0}^{N-1}$ and $\left\{  a_{2,n}%
(x),b_{2,n}(x),\widetilde{a}_{2,n}(x),\widetilde{b}_{2,n}(x)\right\}
_{n=0}^{N-1}$ by substituting the series representations (\ref{phi=}%
)-(\ref{psitil=}) into (\ref{rel 1}), (\ref{rel 2}), (\ref{eig1}) and
(\ref{eig2}). For the first set we obtain%
\begin{gather}
e^{-i\rho_{k}x}\left(  z_{k}+1\right)  \sum_{n=0}^{N-1}\left(  -z_{k}\right)
^{n}b_{1,n}(x)-\mathbf{a}(\rho_{k})e^{-i\rho_{k}x}\left(  \widetilde{z}%
_{k}+1\right)  \sum_{n=0}^{N-1}\left(  -\widetilde{z}_{k}\right)
^{n}\widetilde{a}_{1,n}(x)\nonumber\\
-\mathbf{b}(\rho_{k})e^{i\rho_{k}x}\left(  z_{k}+1\right)  \sum_{n=0}%
^{N-1}\left(  -z_{k}\right)  ^{n}a_{1,n}(x)=\left(  \mathbf{a}(\rho
_{k})-1\right)  e^{-i\rho_{k}x},\quad k=1,\ldots,K, \label{sys11}%
\end{gather}%
\begin{gather}
e^{i\rho_{k}x}\left(  \widetilde{z}_{k}+1\right)  \sum_{n=0}^{N-1}\left(
-\widetilde{z}_{k}\right)  ^{n}\widetilde{b}_{1,n}(x)+\widetilde{\mathbf{a}%
}(\rho_{k})e^{i\rho_{k}x}\left(  z_{k}+1\right)  \sum_{n=0}^{N-1}\left(
-z_{k}\right)  ^{n}a_{1,n}(x)\nonumber\\
-\widetilde{\mathbf{b}}(\rho_{k})e^{-i\rho_{k}x}\left(  \widetilde{z}%
_{k}+1\right)  \sum_{n=0}^{N-1}\left(  -\widetilde{z}_{k}\right)
^{n}\widetilde{a}_{1,n}(x)=\widetilde{\mathbf{b}}(\rho_{k})e^{-i\rho_{k}%
x},\quad k=1,\ldots,K, \label{sys12}%
\end{gather}%
\begin{equation}
e^{-i\rho_{m}x}\left(  z_{m}+1\right)  \sum_{n=0}^{N-1}\left(  -z_{m}\right)
^{n}b_{1,n}(x)-c_{m}e^{i\rho_{m}x}\left(  z_{m}+1\right)  \sum_{n=0}%
^{N-1}\left(  -z_{m}\right)  ^{n}a_{1,n}(x)=-e^{-i\rho_{m}x}, \label{sys13}%
\end{equation}
$m=1,\ldots,M,$%
\begin{equation}
e^{i\widetilde{\rho}_{m}x}\left(  \widetilde{z}_{m}+1\right)  \sum_{n=0}%
^{N-1}\left(  -\widetilde{z}_{m}\right)  ^{n}\widetilde{b}_{1,n}%
(x)-\widetilde{c}_{m}e^{-i\widetilde{\rho}_{m}x}\left(  \widetilde{z}%
_{m}+1\right)  \sum_{n=0}^{N-1}\left(  -\widetilde{z}_{m}\right)
^{n}\widetilde{a}_{1,n}(x)=\widetilde{c}_{m}e^{-i\widetilde{\rho}_{m}x},
\label{sys14}%
\end{equation}
$m=1,\ldots,\widetilde{M}$. Analogously, for the second set of the
coefficients we obtain%
\begin{gather}
e^{-i\rho_{k}x}\left(  z_{k}+1\right)  \sum_{n=0}^{N-1}\left(  -z_{k}\right)
^{n}b_{2,n}(x)-\mathbf{a}(\rho_{k})e^{-i\rho_{k}x}\left(  \widetilde{z}%
_{k}+1\right)  \sum_{n=0}^{N-1}\left(  -\widetilde{z}_{k}\right)
^{n}\widetilde{a}_{2,n}(x)\nonumber\\
-\mathbf{b}(\rho_{k})e^{i\rho_{k}x}\left(  z_{k}+1\right)  \sum_{n=0}%
^{N-1}\left(  -z_{k}\right)  ^{n}a_{2,n}(x)=\mathbf{b}(\rho_{k})e^{i\rho_{k}%
x},\quad k=1,\ldots,K, \label{sys21}%
\end{gather}%
\begin{gather}
e^{i\rho_{k}x}\left(  \widetilde{z}_{k}+1\right)  \sum_{n=0}^{N-1}\left(
-\widetilde{z}_{k}\right)  ^{n}\widetilde{b}_{2,n}(x)+\widetilde{\mathbf{a}%
}(\rho_{k})e^{i\rho_{k}x}\left(  z_{k}+1\right)  \sum_{n=0}^{N-1}\left(
-z_{k}\right)  ^{n}a_{2,n}(x)\nonumber\\
-\widetilde{\mathbf{b}}(\rho_{k})e^{-i\rho_{k}x}\left(  \widetilde{z}%
_{k}+1\right)  \sum_{n=0}^{N-1}\left(  -\widetilde{z}_{k}\right)
^{n}\widetilde{a}_{2,n}(x)=\left(  1-\widetilde{\mathbf{a}}(\rho_{k})\right)
e^{i\rho_{k}x},\quad k=1,\ldots,K, \label{sys22}%
\end{gather}%
\begin{equation}
e^{-i\rho_{m}x}\left(  z_{m}+1\right)  \sum_{n=0}^{N-1}\left(  -z_{m}\right)
^{n}b_{2,n}(x)-c_{m}e^{i\rho_{m}x}\left(  z_{m}+1\right)  \sum_{n=0}%
^{N-1}\left(  -z_{m}\right)  ^{n}a_{2,n}(x)=c_{m}e^{i\rho_{m}x}, \label{sys23}%
\end{equation}
$m=1,\ldots,M,$%
\begin{equation}
e^{i\widetilde{\rho}_{m}x}\left(  \widetilde{z}_{m}+1\right)  \sum_{n=0}%
^{N-1}\left(  -\widetilde{z}_{m}\right)  ^{n}\widetilde{b}_{2,n}%
(x)-\widetilde{c}_{m}e^{-i\widetilde{\rho}_{m}x}\left(  \widetilde{z}%
_{m}+1\right)  \sum_{n=0}^{N-1}\left(  -\widetilde{z}_{m}\right)
^{n}\widetilde{a}_{2,n}(x)=e^{i\widetilde{\rho}_{m}x}, \label{sys24}%
\end{equation}
$m=1,\ldots,\widetilde{M}$.

Observe that the matrices of the systems (\ref{sys11})-(\ref{sys14}) and
(\ref{sys21})-(\ref{sys24}) coincide, while the vectors on the right-hand side
are different. Thus, for any $x\in\mathbb{R}$ we have two systems of linear
algebraic equations
\begin{equation}
AX_{1}=B_{1}\quad\text{and}\quad AX_{2}=B_{2}, \label{systems}%
\end{equation}
where
\begin{equation}
X_{1}=\left(  b_{1,0}(x),\ldots,b_{1,N-1}(x),\widetilde{b}_{1,0}%
(x),\ldots,\widetilde{b}_{1,N-1}(x),a_{1,0}(x),\ldots,a_{1,N-1}(x),\widetilde
{a}_{1,0}(x),\ldots,\widetilde{a}_{1,N-1}(x)\right)  ^{T}, \label{X1}%
\end{equation}%
\begin{equation}
X_{2}=\left(  b_{2,0}(x),\ldots,b_{2,N-1}(x),\widetilde{b}_{2,0}%
(x),\ldots,\widetilde{b}_{2,N-1}(x),a_{2,0}(x),\ldots,a_{2,N-1}(x),\widetilde
{a}_{2,0}(x),\ldots,\widetilde{a}_{2,N-1}(x)\right)  ^{T}, \label{X2}%
\end{equation}
and the matrix $A$ has the following form. For $k=1,\ldots,K$,%
\[
A_{kj}=\left(  -1\right)  ^{j-1}e^{-i\rho_{k}x}\left(  z_{k}+1\right)
z_{k}^{j-1},\,j=1,\ldots,N,\quad\quad A_{kj}=0,\,j=N+1,\ldots,2N,
\]%
\[
A_{kj}=-\mathbf{b}(\rho_{k})e^{i\rho_{k}x}\left(  z_{k}+1\right)  \left(
-z_{k}\right)  ^{j-(2N+1)},\,j=2N+1,\ldots,3N,
\]%
\[
A_{kj}=-\mathbf{a}(\rho_{k})e^{-i\rho_{k}x}\left(  \widetilde{z}_{k}+1\right)
\left(  -\widetilde{z}_{k}\right)  ^{j-(3N+1)},\,j=3N+1,\ldots,4N.
\]
For $k=K+1,\ldots,2K$,%
\[
A_{kj}=0,\,j=1,\ldots,N,\quad\quad A_{kj}=e^{i\rho_{k-K}x}\left(
\widetilde{z}_{k-K}+1\right)  \left(  -\widetilde{z}_{k-K}\right)
^{j-(N+1)},\,j=N+1,\ldots,2N,
\]%
\[
A_{kj}=\widetilde{\mathbf{a}}(\rho_{k-K})e^{i\rho_{k-K}x}\left(
z_{k-K}+1\right)  \left(  -z_{k-K}\right)  ^{j-(2N+1)},\,j=2N+1,\ldots,3N,
\]%
\[
A_{kj}=-\widetilde{\mathbf{b}}(\rho_{k-K})e^{-i\rho_{k-K}x}\left(
\widetilde{z}_{k-K}+1\right)  \left(  -\widetilde{z}_{k-K}\right)
^{j-(3N+1)},\,j=3N+1,\ldots,4N.
\]
For $k=2K+1,\ldots,2K+M$ (identifying $\rho_{m}$, $z_{m}$ and $c_{m}$ for
$m=1,\ldots,M$ with $\rho_{k}$, $z_{k}$ and $c_{k}$ for $k=2K+1,\ldots,2K+M$)
we have
\[
A_{kj}=e^{-i\rho_{k}x}\left(  z_{k}+1\right)  \left(  -z_{k}\right)
^{j-1},\,j=1,\ldots,N,\quad\quad A_{kj}=0,\,j=N+1,\ldots,2N,
\]%
\[
A_{kj}=-c_{k}e^{i\rho_{k}x}\left(  z_{k}+1\right)  \left(  -z_{k}\right)
^{j-(2N+1)},\,j=2N+1,\ldots,3N,\quad\quad A_{kj}=0,\,j=3N+1,\ldots,4N.
\]
For $k=2K+M+1,\ldots,2K+M+\widetilde{M}$ (identifying $\widetilde{\rho}_{m}$,
$\widetilde{z}_{m}$ and $\widetilde{c}_{m}$ for $m=1,\ldots,\widetilde{M}$
with $\rho_{k}$, $z_{k}$ and $c_{k}$ for $k=2K+M+1,\ldots,2K+M+\widetilde{M}$)
we have%
\[
A_{kj}=0,\,j=1,\ldots,N,\quad\quad A_{kj}=e^{i\rho_{k}x}\left(  z_{k}%
+1\right)  \left(  -z_{k}\right)  ^{j-(N+1)},\,j=N+1,\ldots,2N,
\]%
\[
A_{kj}=0,\,j=2N+1,\ldots,3N,\quad\quad A_{kj}=-c_{k}e^{-i\rho_{k}x}\left(
z_{k}+1\right)  \left(  -z_{k}\right)  ^{j-(3N+1)},\,j=3N+1,\ldots,4N.
\]
The right-hand side vectors have the form%
\[
B_{1,k}=\left(  \mathbf{a}(\rho_{k})-1\right)  e^{-i\rho_{k}x},\,k=1,\ldots
,K,\quad B_{1,k}=\widetilde{\mathbf{b}}(\rho_{k-K})e^{-i\rho_{k-K}%
x},\,k=K+1,\ldots,2K,
\]%
\[
B_{1,k}=-e^{-i\rho_{k}x},\,k=2K+1,\ldots,2K+M,\quad B_{1,k}=c_{k}e^{-i\rho
_{k}x},\,k=2K+M+1,\ldots,2K+M+\widetilde{M}%
\]
and%
\[
B_{2,k}=\mathbf{b}(\rho_{k})e^{i\rho_{k}x},\,k=1,\ldots,K,\quad B_{2,k}%
=\left(  1-\widetilde{\mathbf{a}}(\rho_{k-K})\right)  e^{i\rho_{k-K}%
x},\,k=K+1,\ldots,2K,
\]%
\[
B_{2,k}=c_{k}e^{i\rho_{k}x},\,k=2K+1,\ldots,2K+M,\quad B_{2,k}=e^{i\rho_{k}%
x},\,k=2K+M+1,\ldots,2K+M+\widetilde{M}.
\]

Solving systems (\ref{systems}) aims at finding any of the first SPPS
coefficients, for example, $b_{0}(x)$. As we explained at the beginning of
this section, having computed it at a set of points allows us to recover
$q(x)$ and $r(x)$ from (\ref{q and r}).

\section{Numerical examples}

The proposed approach can be implemented directly using an available numerical
computing environment. All the reported computations were performed in Matlab
R2024a on an Intel i7-1360P equipped laptop computer and took no more than
several seconds. For all the integrations involved we used the Newton--Cottes
six point integration rule with $2500$ nodes per unit. Numerical
differentiation which is involved in the last step for recovering $q(x)$ and
$r(x)$ from the first SPPS coefficients ($r(x)$ by the second formula in
(\ref{q and r}) and $q(x)$ by (\ref{q from phitil})) was done numerically with
the aid of spline interpolation (Matlab routine `\emph{spapi}') and
differentiation of the spline (Matlab routine `\emph{fnder}').

\subsection{Direct scattering}

\textbf{Example 1. }We begin with the Zakharov-Shabat system with a potential,
for which the analytical expressions of the scattering data are known
\cite{Tovbis et al 2004}, \cite{Trogdon2021}%
\[
q(x)=-iA\operatorname*{sech}(x)\exp(-i\gamma A\ln\cosh(x)),\quad
A>0,\,\gamma\in\mathbb{R},
\]%
\[
r(x)=-\overline{q(x)}.
\]
System (\ref{ZS1}), (\ref{ZS2}) with these potentials gives rise the
scattering data
\[
\mathbf{a}(\rho)=\frac{\Gamma\left(  \omega\left(  \rho\right)  \right)
\,\Gamma\left(  \omega\left(  \rho\right)  -\omega_{-}-\omega_{+}\right)
}{\Gamma\left(  \omega\left(  \rho\right)  -\omega_{+}\right)  \,\Gamma\left(
\omega\left(  \rho\right)  -\omega_{-}\right)  },
\]%
\[
\mathbf{b}(\rho)=iA^{-1}2^{-i\gamma A}\frac{\Gamma\left(  \omega\left(
\rho\right)  \right)  \,\Gamma\left(  1-\omega\left(  \rho\right)  +\omega
_{-}+\omega_{+}\right)  }{\Gamma\left(  \omega_{+}\right)  \,\Gamma\left(
\omega_{-}\right)  },
\]%
\[
\omega\left(  \rho\right)  =-i\rho-A\gamma\frac{i}{2}+\frac{1}{2},\quad
\omega_{+}=-iA\left(  T+\frac{\gamma}{2}\right)  ,
\]%
\[
\omega_{-}=iA\left(  T-\frac{\gamma}{2}\right)  ,\quad T=\sqrt{\frac
{\gamma^{2}}{4}-1}.
\]
The zeros of $\mathbf{a}(\rho)$ are given by
\[
\rho_{m}=AT-i(m-\frac{1}{2}),\quad m=1,2,\ldots,M\text{\quad with
}M=\left\lfloor \frac{1}{2}+A\left\vert T\right\vert \right\rfloor .
\]
Also $c_{m}=\mathbf{b}(\rho_{m})$. For the numerical test, following
\cite{Trogdon2021}, we choose $A=1.65$ and $\gamma=0.1$. To compute the
scattering data, $160$ SPPS coefficients for each Jost solution were computed
following the recurrent integration procedures devised in Section
\ref{Sect recurrent}. In Fig. 1 the decay of the coefficients $a_{1,n}(0)$,
$a_{2,n}(0)$ is presented.%

\begin{figure}
[ptbh]
\begin{center}
\includegraphics[
height=3.6625in,
width=4.8758in
]%
{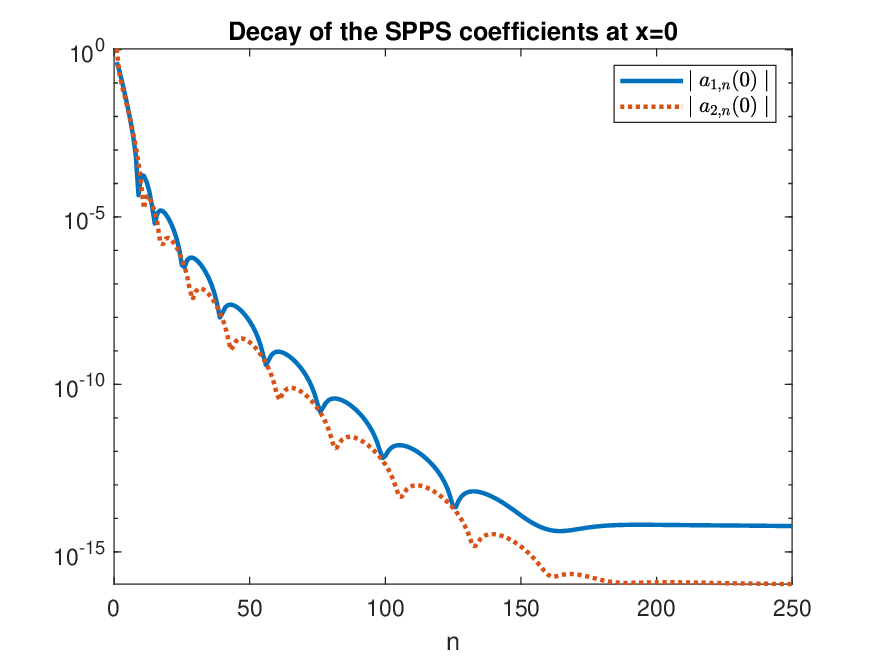}%
\caption{Decay of the computed coefficients $a_{1,n}(0)$ and $a_{2,n}(0)$
corresponding to Example 1, when $n$ increases.}%
\label{FigDecayEx1}%
\end{center}
\end{figure}

In Fig. 2 we compare the computed functions $\mathbf{a}(\rho)$ and
$\mathbf{b}(\rho)$ with their analytical expressions. The resulting absolute
error is of order $10^{-13}$. This is at least two orders more accurate than
the reported most accurate results in \cite{Trogdon2021}. We emphasize that
the SPPS series representations allow us to compute the matrix $S(\rho)$ on
immense intervals of $\rho\in\mathbb{R}$ and in no time, with a
non-deteriorating accuracy and at an arbitrarily large number of points. This
is because having computed once the SPPS coefficients at $x=0$, all further
computations reduce to computing values of polynomials in $z$ (and
$\widetilde{z}$ in examples below) when $z$ (and $\widetilde{z}$) runs along
the unit circle (which corresponds to $\rho$ running along the real axis).%

\begin{figure}
[ptbh]
\begin{center}
\includegraphics[
height=4.2964in,
width=5.7199in
]%
{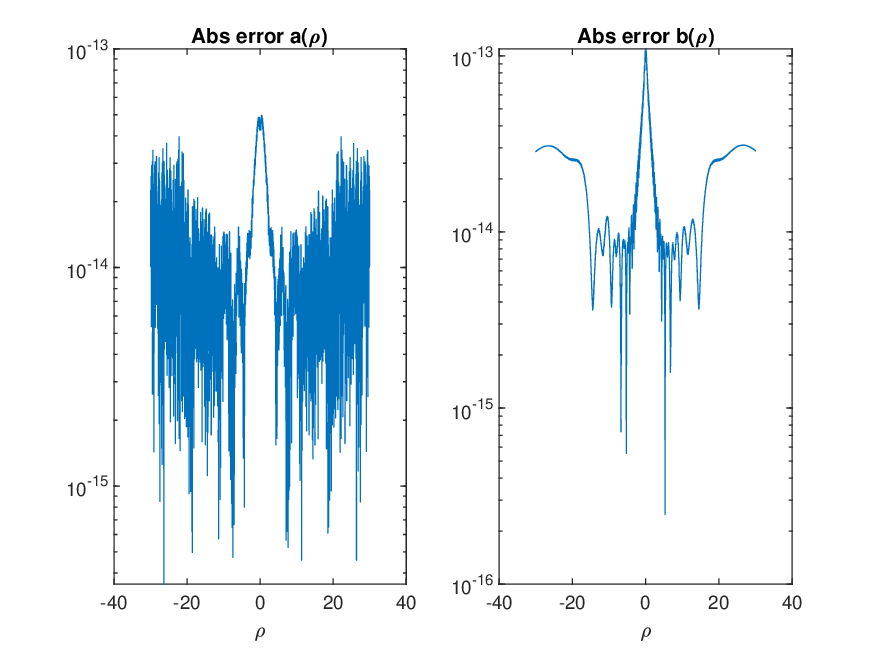}%
\caption{Absolute error of the computed $\mathbf{a}(\rho)$ (left) and
$\mathbf{b}(\rho)$ (right). The maximum absolute error of $\mathbf{a}(\rho)$
resulted in $4.9\cdot10^{-14}$, and that of $\mathbf{b}(\rho)$ in
$1.09\cdot10^{-13}$.}%
\label{Fig2}%
\end{center}
\end{figure}

The two eigenvalues $\rho_{1}\approx0.14793620932365i$ and $\rho_{2}%
\approx1.14793620932365i$ were computed with the absolute error $1.02\cdot
10^{-14}$ and $4.06\cdot10^{-14}$, respectively. We mention that for this, the
location of the roots of the polynomial of the variable $z$ (see Section
\ref{SectDirect}) was performed with the aid of the Matlab routine
`\emph{roots}'. We emphasize that neither in this example nor in the
subsequent ones do spurious eigenvalues occur. The corresponding norming
constants $c_{1}\approx-0.187821133726638+0.982203248684122i$ and
$c_{2}\approx-0.0643040290406992-0.997930354207713i$ were computed with the
absolute error $7.9\cdot10^{-14}$ and $2.44\cdot10^{-14}$, respectively.\ 

\bigskip

\textbf{Example 2. }Consider the case \cite{Trogdon2021}
\begin{equation}
q(x)=e^{-x^{2}},\quad r(x)=-2e^{-x^{2}+ix}. \label{Potential Trogdon}%
\end{equation}
The analytical expressions for the corresponding scattering data are unknown,
however we can compare the results of our computation with those given in
\cite{Trogdon2021}. Moreover, the attained accuracy becomes more evident when
we solve the inverse problem using the obtained results as the input data for
the inverse problem (see Example 4 below).

To compute the scattering data four hundreed SPPS coefficients for each of the
Jost solutions were computed following the recurrent integration procedures
devised in Section \ref{Sect recurrent}. In Fig. 3 the decay of the
coefficients $a_{1,n}(0)$, $a_{2,n}(0)$ is presented.%

\begin{figure}
[ptbh]
\begin{center}
\includegraphics[
height=4.1321in,
width=5.5011in
]%
{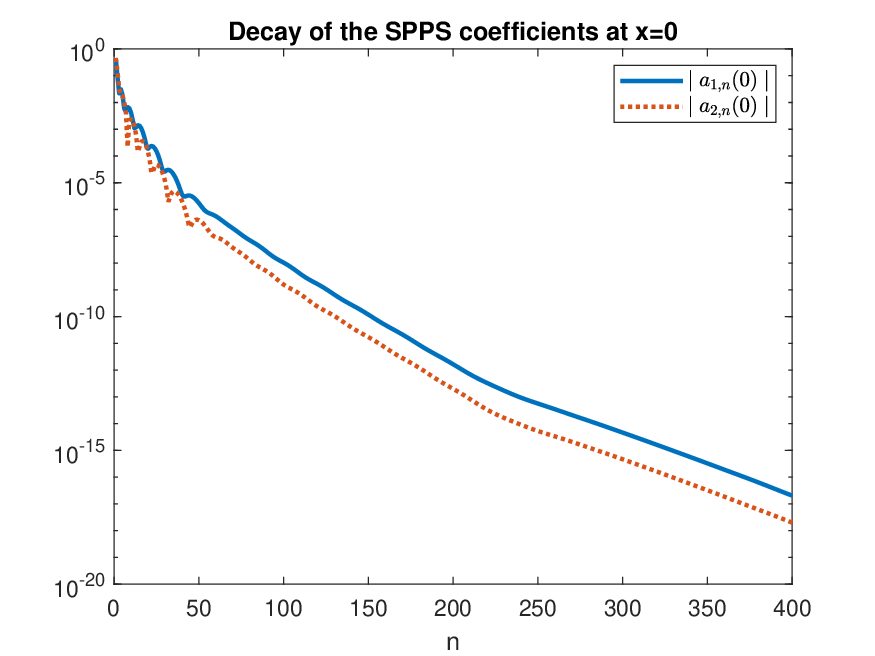}%
\caption{Decay of the computed coefficients $a_{1,n}(0)$ and $a_{2,n}(0)$
corresponding to Example 2, when $n$ increases.}%
\label{FigDecayEx2}%
\end{center}
\end{figure}

The values of $\mathbf{a}(\rho_{k})$, $\widetilde{\mathbf{a}}(\rho_{k})$,
$\mathbf{b}(\rho_{k})$ and $\widetilde{\mathbf{b}}(\rho_{k})$ are computed at
$4000$ points $\rho_{k}$ distributed uniformly along the segment $[-30,30]$.
Equality (\ref{aatil}) may serve as an additional reliable indicator of the
accuracy of the computed entries of the matrix $S(\rho)$. In Fig. 4 we show
that (\ref{aatil}) is fulfilled indeed with a remarkable accuracy (here the
maximum value of the depicted expression is $1.2\cdot10^{-14}$).%

\begin{figure}
[ptbh]
\begin{center}
\includegraphics[
height=4.3785in,
width=5.8288in
]%
{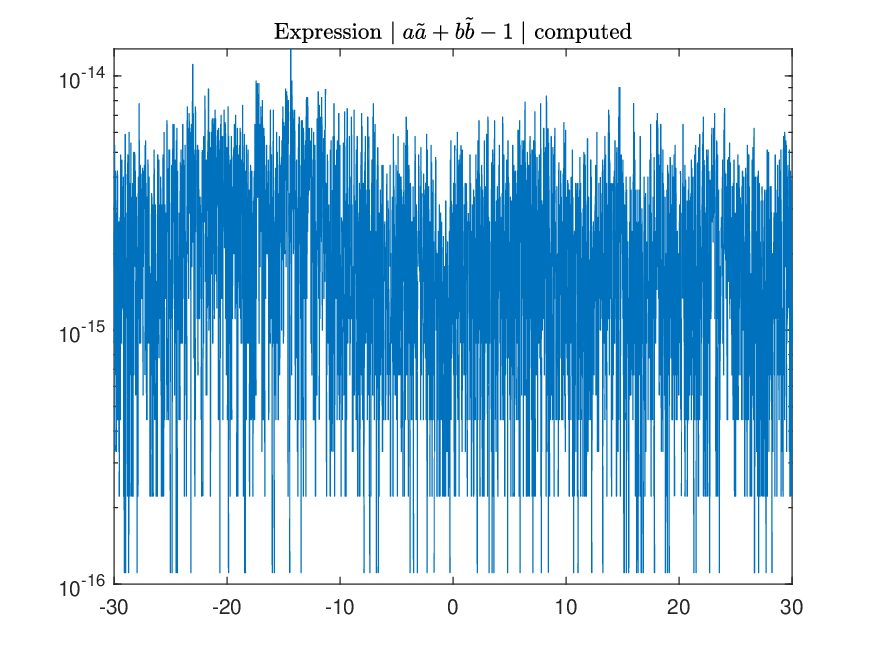}%
\caption{Equality (\ref{aatil}) is used as an indicator of the accuracy of the
computed entries of the matrix $S(\rho)$.}%
\label{FigSquaresEx2}%
\end{center}
\end{figure}

Two eigenvalues are detected%
\begin{equation}
\rho_{1}\approx0.250000000000056+0.501700389937887i, \label{rho1}%
\end{equation}%
\begin{equation}
\widetilde{\rho}_{1}\approx0.250000000000051-0.501700389937864i.
\label{rhotil1}%
\end{equation}
In \cite{Trogdon2021} also two eigenvalues are reported%
\[
\rho_{1}\approx0.25+0.517003899379i,\quad\widetilde{\rho}_{1}\approx
0.25+0.517003899379i.
\]
Since the digits in the imaginary parts coincide with (\ref{rho1}),
(\ref{rhotil1}) except for the missing zero in the hundredths digit place we
attribute this discrepancy to a typo in \cite{Trogdon2021}.

\textbf{Example 3. }Consider the case
\begin{equation}
q(x)=\pi\exp(-x^{2}+i\sin(\pi x)), \label{Pot3q}%
\end{equation}%
\begin{equation}
r(x)=-\pi\exp(-x^{2}-i\cos(\pi x)). \label{Pot3r}%
\end{equation}
Analytical expressions for the corresponding scattering data are unknown.

To compute the scattering data seven hundreed SPPS coefficients for each of
the Jost solutions were computed following the recurrent integration
procedures from Section \ref{Sect recurrent}. In Fig. 5 the decay of the
coefficients $a_{1,n}(0)$, $a_{2,n}(0)$ is presented.%

\begin{figure}
[ptb]
\begin{center}
\includegraphics[
height=3.7291in,
width=4.964in
]%
{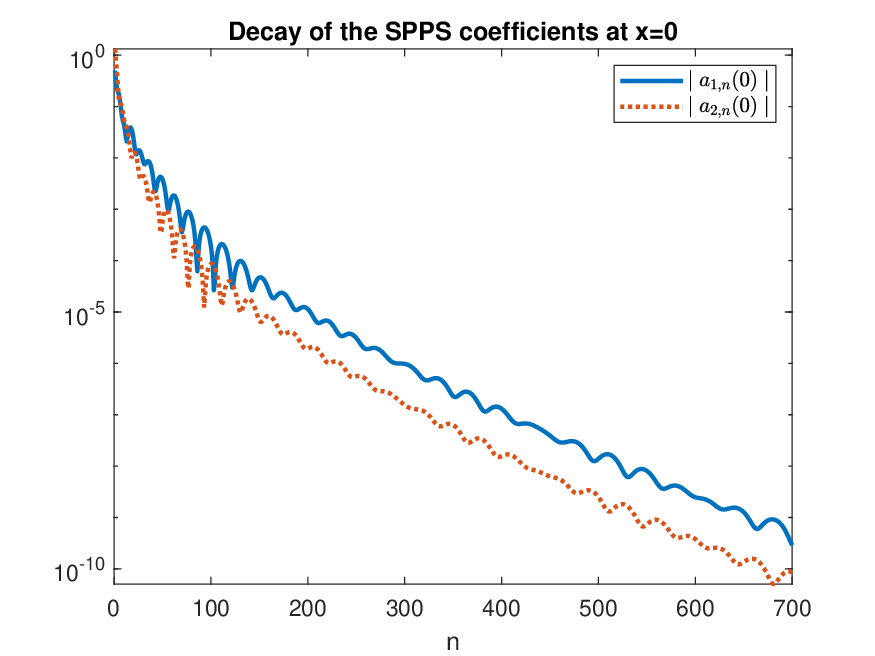}%
\caption{Decay of the computed coefficients $a_{1,n}(0)$ and $a_{2,n}(0)$
corresponding to Example 3, when $n$ increases.}%
\label{FigDecayEx3New}%
\end{center}
\end{figure}

The values of $\mathbf{a}(\rho_{k})$, $\widetilde{\mathbf{a}}(\rho_{k})$,
$\mathbf{b}(\rho_{k})$ and $\widetilde{\mathbf{b}}(\rho_{k})$ are computed at
$5000$ points $\rho_{k}$ chosen as follows. Half of them are taken in the form
$\rho_{k}=10^{\alpha_{k}}$ with $\alpha_{k}$ being distributed uniformly on
$\left[  \log(0.001),\log(70)\right]  \approx\left[  -3,\,1.845\right]  $.
This is a logarithmically spaced point distribution on the segment
$[0.001,70]$, with points more densely distributed near $\rho=0.001$ and less
densely distributed as $\rho$ increases. The other half is taken symmetrically
with respect to zero. This logarithmically spaced point distribution for the
scattering data computed may lead to more accurate results when solving the
inverse problem, see Example 6 below. Again we use equality (\ref{aatil}) for
checking the accuracy of the computed entries of the matrix $S(\rho)$, see
Fig. 6. We observe that the entries of $S(\rho)$ are computed with a lower
accuracy as compared to the previous example, which can be explained by a
lower decay of the SPPS coefficients, see Fig. 5.%

\begin{figure}
[ptbh]
\begin{center}
\includegraphics[
height=3.9557in,
width=5.2667in
]%
{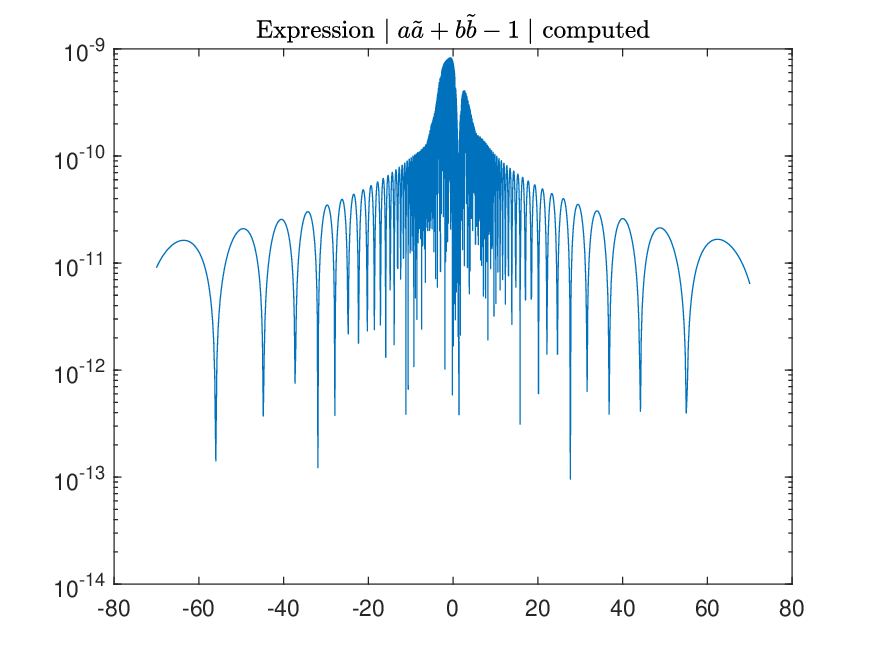}%
\caption{Equality (\ref{aatil}) is used as an indicator of the accuracy of the
computed entries of the matrix $S(\rho)$ corresponding to Example 3.}%
\label{FigSquaresEx3}%
\end{center}
\end{figure}

Four eigenvalues are detected%
\[
\rho_{1}\approx0.281405857470267+1.94356920198665i,\quad\rho_{2}%
\approx0.545035754764913+0.51356582669352i,
\]
\[
\widetilde{\rho}_{1}\approx1.18535887013205-0.0492419359968676i,\quad
\widetilde{\rho}_{2}\approx-1.98047598318108-0.87978108360884i
\]
with the corresponding norming constants
\begin{align*}
c_{1}  &  \approx-1.02697836088912-0.501965022137801i,\\
c_{2}  &  \approx0.482038890052823+0.306291387212402i,\\
\widetilde{c}_{1}  &  \approx1.67939149873166-0.553134958681549i,\\
\widetilde{c}_{2}  &  \approx-0.156577573656642-0.667640588064389i.
\end{align*}

In Example 6 below the computed scattering data are used for recovering the
potentials (\ref{Pot3q}) and (\ref{Pot3r}).

\subsection{Inverse scattering}

Solution of the inverse scattering problems consists in solving two systems
(\ref{systems}) at a sufficiently large number of points $x_{j}\in
\lbrack-l,l]$ for some $l>0$ and recovering the potentials $q(x)$ and $r(x)$
from the first SPPS coefficients, as explained in the beginning of Section
\ref{Sect Inverse} ($q(x)$ from (\ref{q from phitil}) and $r(x)$ from the
second equality in (\ref{q and r})). In the first two examples we choose $l=5$
and $N=50$ in (\ref{X1}), (\ref{X2}).

\textbf{Example 4. }We use the computed scattering data from Example 1 as
input data for the inverse scattering problem. First we take the values of
$\mathbf{a}(\rho_{k})$ and $\mathbf{b}(\rho_{k})$ computed at $4000$ points
$\rho_{k}$ distributed uniformly along the segment $[-30,30]$, complemented by
the computed discrete scattering data. In Fig. 7 the recovered potential is presented.%

\begin{figure}
[ptbh]
\begin{center}
\includegraphics[
height=4.3128in,
width=5.7406in
]%
{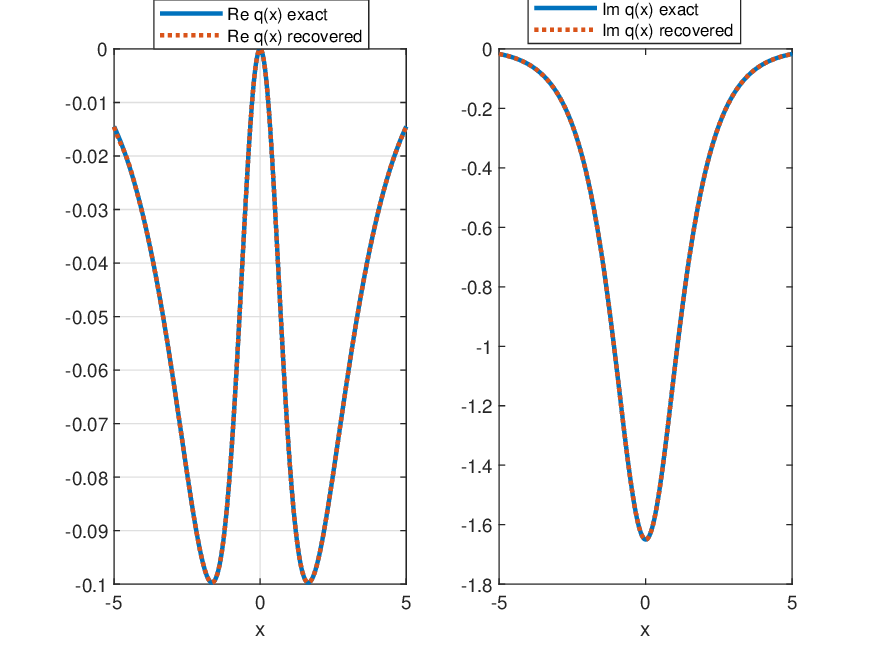}%
\caption{The potential from Example 4 is recovered from the scattering data
computed in Example 1. The maximum absolute error of the recovered potential
resulted in $7.9\cdot10^{-5}$. }%
\label{FigPot1}%
\end{center}
\end{figure}
The  maximum absolute error of the recovered potential resulted in
$7.9\cdot10^{-5}$. 

When more data are available the recovery is more accurate. As an example we
computed the input data at $14\cdot10^{3}$ points $\rho_{k}$ distributed
uniformly along the segment $[-130,130]$. The maximum absolute error of the
recovered potential resulted in $7.8\cdot10^{-6}$.

\textbf{Example 5. }We use the computed scattering data from Example 2 as
input data for the inverse scattering problem. First we take the values of
$\mathbf{a}(\rho_{k})$ and $\mathbf{b}(\rho_{k})$ computed at $4000$ points
$\rho_{k}$ distributed uniformly along the segment $[-30,30]$, complemented by
the computed discrete scattering data. In Figs. 8 and 9 the recovered
potentials are presented.%

\begin{figure}
[ptbh]
\begin{center}
\includegraphics[
height=4.1364in,
width=5.5063in
]%
{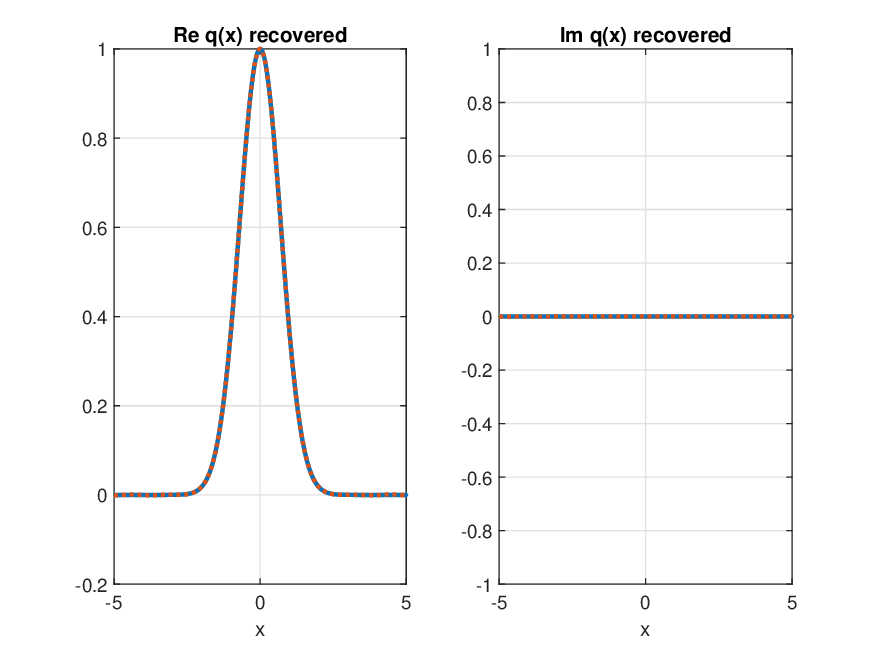}%
\caption{Recovered potential $q(x)$ from Example 5. The maximum absolute error
resulted in $1.05\cdot10^{-3}$}%
\label{FigPot2q}%
\end{center}
\end{figure}
%

\begin{figure}
[ptbh]
\begin{center}
\includegraphics[
height=3.9946in,
width=5.3186in
]%
{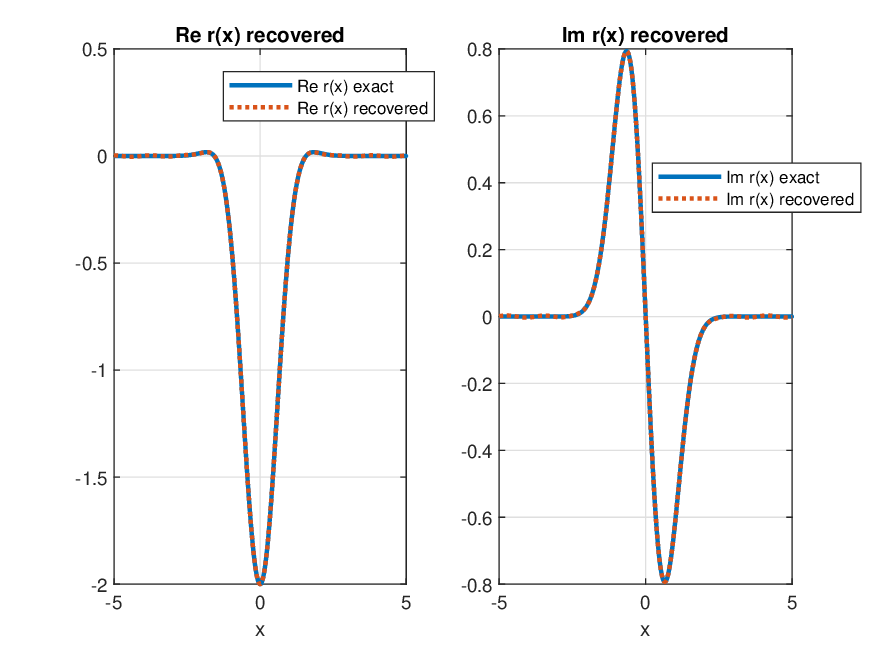}%
\caption{Recovered potential $r(x)$ from Example 5. The maximum absolute error
resulted in $3.39\cdot10^{-3}$.}%
\label{FigPot2r}%
\end{center}
\end{figure}

When more data are available the recovery is more accurate. As an example we
computed the input data at $14\cdot10^{3}$ points $\rho_{k}$ distributed
uniformly along the segment $[-130,130]$. The maximum absolute error of the
recovered potential resulted in $7.9\cdot10^{-4}$.

\textbf{Example 6. }We use the computed scattering data from Example 3 to
reconstruct the potentials (\ref{Pot3q}), (\ref{Pot3r}), see Figs. 10, 11. The
maximum absolute error of the reconstruction of $q(x)$ and \ $r(x)$ resulted
in $0.037$ and $0.045$, respectively. The parameter $N$ in (\ref{X1}),
(\ref{X2}) was chosen as $N=90$.%

\begin{figure}
[ptb]
\begin{center}
\includegraphics[
height=3.9133in,
width=5.2096in
]%
{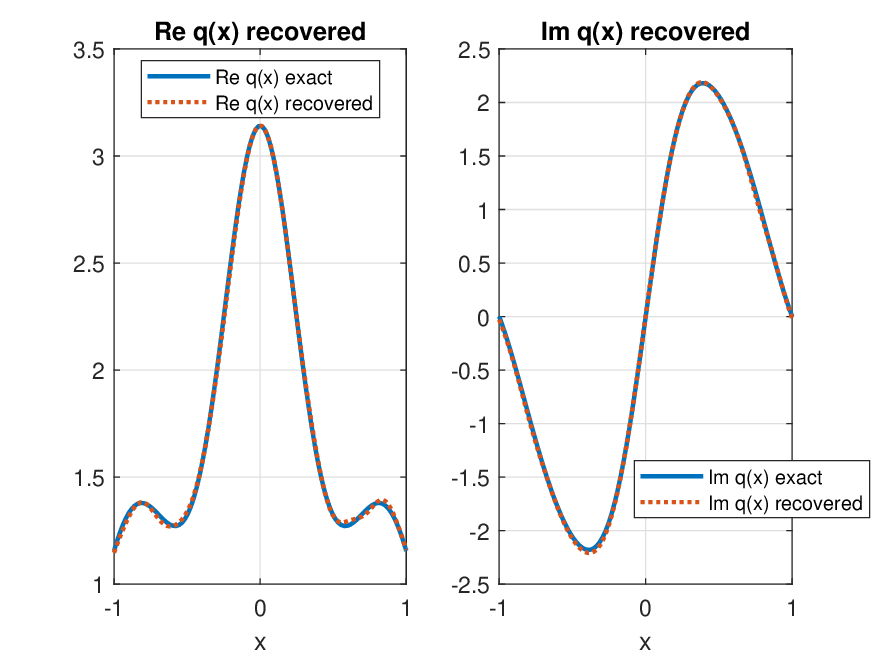}%
\caption{Recovered potential $q(x)$ from Example 6. The maximum absolute error
resulted in  $0.037$}%
\label{FigPot3q}%
\end{center}
\end{figure}

%

\begin{figure}
[ptb]
\begin{center}
\includegraphics[
height=3.8207in,
width=5.0842in
]%
{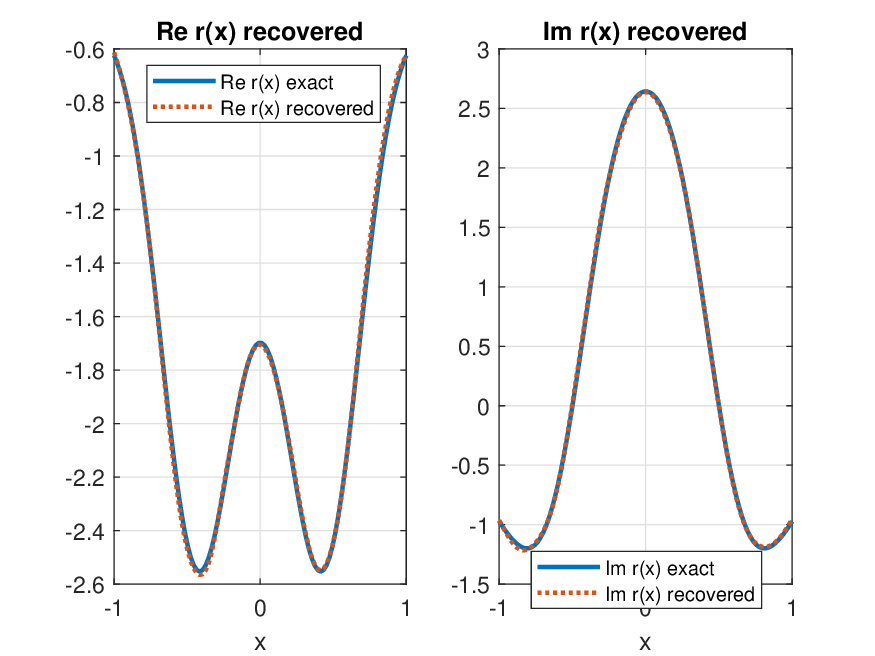}%
\caption{Recovered potential $r(x)$ from Example 6. The maximum absolute error
resulted in  $0.045$,}%
\label{FigPot3r}%
\end{center}
\end{figure}

The result of the computation shows that the solution of both the direct and
inverse problems is accurate. When the analytical expression of the potentials
is unknown, useful indicators of the accuracy\ of the computed SPPS
coefficients can be easily derived, for example, from (\ref{Wronski}).

\section{Conclusions\label{Sect Conclusions}}

The AKNS system with complex valued, sufficiently fast decaying potentials is
considered. For its Jost solutions power series representations are obtained
in terms of a mapped spectral parameter. The solution of the direct scattering
problem thus reduces to computing the coefficients of the series and locating
zeros of resulting polynomials inside of the unit disk. In its turn the
inverse scattering problem reduces to the solution of two systems of linear
algebraic equations, one for the first components of the Jost solutions and
the other for the second ones. The potentials are recovered from the first
elements of the solution vectors. The proposed method for solving both the
direct and inverse problems leads to easily implemented, direct and accurate
algorithms. Their application to the solution of initial value problems for
the nonlinear Schr\"{o}dinger equation is developed in \cite{Kr NLS arxiv}.

\textbf{Funding information }Ministry of Education and Science of Russia,
under Grant Agreement No. 075-02-2025-1720.

\textbf{Data availability} The data that support the findings of this study
are available upon reasonable request.

\textbf{Conflict of interest }This work does not have any conflict of interest.

\end{document}